\documentclass[12pt]{article}
\usepackage{amssymb}
\usepackage{amsmath}
\usepackage{hyperref}

\def\sqr#1#2{{\vcenter{\vbox{\hrule height.#2pt
              \hbox{\vrule width.#2pt height#1pt \kern#1pt \vrule width.#2pt}
              \hrule height.#2pt}}}}

\def\dbR{{\mathop{\rm l\negthinspace R}}}
\def\dbC{{\mathop{\rm l\negthinspace\negthinspace\negthinspace C}}}

\def \Ker{\mbox{\rm Ker}}
\def\3n{\negthinspace \negthinspace \negthinspace }
\def\2n{\negthinspace \negthinspace }
\def\1n{\negthinspace }

\def\dbC{{\mathbb{C}}}

\def\dbR{{\mathbb{R}}}

\def\dbZ{{\mathbb{Z}}}
\def\={\buildrel \triangle \over =}

\def\ov{\overline  }
\def\a{\alpha}
\def\b{\beta}
\def\g{\gamma}
\def\d{\delta}

\def\k{\kappa}
\def\l{\lambda}

\def\si{\sigma}

\def\f{\varphi}
\def\th{\theta}
\def\o{\omega}

\def\ns{\noalign{\ss} }

\def\pa{\partial}
\def\R{\mathbb{R}}

\def\cB{{\cal B}}

\def\cE{{\cal E}}
\def\cF{{\cal F}}
\def\cG{{\cal G}}

\def\cJ{{\cal J}}
\def\cK{{\cal K}}
\def\cL{{\cal L}}

\def\cN{{\cal N}}

\def\no{\noindent}

\def\ms{\medskip}
\def\bs{\bigskip}
\def\q{\quad}
\def\qq{\qquad}

\def\max{\mathop{\rm max}}
\def\min{\mathop{\rm min}}

\def\pa{\partial}

\def\wt{\widetilde}
\def\cd{\cdot}
\def\cds{\cdots}

\def\spa{\hbox{\rm span$\,$}}

\def\wh{\widehat}

\def\|{\Big |}
\def\({\Big (}
\def\){\Big )}
\def\[{\Big[}
\def\]{\Big]}
\def\be{\begin{equation}}
\def\bel{\begin{equation}\label}
\def\ee{\end{equation}}
\def\bt{\begin{theorem}}
\def\bcd{\begin{condition}}
\def\ecd{\end{condition}}
\def\et{\end{theorem}}
\def\bc{\begin{corollary}}
\def\ec{\end{corollary}}
\def\bde{\begin{definition}}
\def\ede{\end{definition}}
\def\bl{\begin{lemma}}
\def\el{\end{lemma}}
\def\bp{\begin{proposition}}
\def\ep{\end{proposition}}
\def\br{\begin{remark}}
\def\er{\end{remark}}
\def\ba{\begin{array}}
\def\ea{\end{array}}
\def\ed{\end{document}}
\def\ns{\noalign{\ms}}
\def\ds{\displaystyle}

\def\square#1{\vbox{\hrule\hbox{\vrule height#1%
     \kern#1\vrule}\hrule}}
\def\rectangle#1#2{\vbox{\hrule\hbox{\vrule height#1%
     \kern#2\vrule}\hrule}}

\font\tenbb=msbm10 \font\sevenbb=msbm7
\font\fivebb=msbm5

\newfam\bbfam
\scriptscriptfont\bbfam=\fivebb
\textfont\bbfam=\tenbb
\scriptfont\bbfam=\sevenbb

\newtheorem{lemma}{Lemma}[section]
\newtheorem{remark}{Remark}[section]

\newtheorem{theorem}{Theorem}[section]
\newtheorem{corollary}{Corollary}[section]

\newtheorem{definition}{Definition}[section]
\newtheorem{proposition}{Proposition}[section]
\newtheorem{condition}{Condition}[section]

\makeatletter
   
   \@addtoreset{equation}{section}
\makeatother

\begin{document}

\title{\bf Local rapid stabilization for a Korteweg-de Vries equation with a Neumann boundary control on the right}

\author{Jean-Michel Coron\thanks{Sorbonne Universit\'{e}s, UPMC Univ Paris 06, UMR 7598, Laboratoire
Jacques-Louis Lions, F-75005, Paris, France
E-mail: \texttt{coron@ann.jussieu.fr}. JMC
was supported by ERC advanced grant 266907
(CPDENL) of the 7th Research Framework
Programme (FP7).} \q and \q Qi
L\"{u}\thanks{School of Mathematical
Sciences, University of Electronic Science
and Technology of China, Chengdu 610054,
China. E-mail: \texttt{luqi59@163.com}. QL
was supported by ERC advanced grant 266907
(CPDENL) of the 7th Research Framework
Programme (FP7), the NSF of China under
grant 11101070 and the Fundamental Research
Funds for the Central Universities in China
under grants ZYGX2012J115.}}

\date{}

\maketitle

\begin{abstract}\no
This paper is devoted to the study of the
rapid exponential stabilization problem for
a controlled Korteweg-de Vries equation on a
bounded interval with homogeneous Dirichlet
boundary conditions and Neumann boundary
control at the right endpoint of the
interval. For every noncritical length, we
build a feedback control law  to force the
solution of the closed-loop system to decay
exponentially to zero with arbitrarily
prescribed decay rates, provided that the
initial datum is small enough. Our approach
relies on the construction of a suitable
integral transform and can be applied to
many other equations.
\end{abstract}

\bs

\no{\bf 2010 Mathematics Subject
Classification}. 93D15, 35Q53. \bs

\no{\bf Key Words}.  Korteweg-de Vries equation,
stabilization,  integral transform.

 \ms
\section{Introduction }


Let $L>0$. Consider the following
controlled Korteweg-de Vries (KdV) equation
on $(0,L)$
\begin{equation}\label{csystem1}
\left\{
\begin{array}{ll}\ds
v_t + v_x + v_{xxx} + vv_x = 0 &\mbox{ in }
(0,+\infty)\times (0,L),\\
\ns\ds v(t,0) = v(t,L)=0 &\mbox{ on }
(0,+\infty),\\
\ns\ds v_x(t,L)=f(t)&\mbox{ on }
(0,+\infty).
\end{array}
\right.
\end{equation}
This is a control system, where, at time
$t\in [0,+\infty)$, the state is
$v(t,\cd)\in L^2(0,L)$ and  the control is
$f(t)\in \R$.

We are concerned with the following
stabilization problem for the system
\eqref{csystem1}.

\vspace{0.3cm}

Problem {\bf (S)}: Let $\l>0$. Does there
exist a linear feedback control $F:
L^2(0,L) \to \mathbb{R}$ such that, for
some $\d>0$, every solution $v$ of
\eqref{csystem1} with $f(t)=F(v(t,\cd))$
satisfies
\begin{equation}
\label{est-sur-v}
 |v(t,\cd)|_{L^2(0,L)} \leq
Ce^{-\l t}|v(0,\cd)|_{L^2(0,L)},
\end{equation}
for some $C>0$, provided that
$|v(0,\cd)|_{L^2(0,L)}\leq \d$?

Let us point out that, for every continuous
linear map $F:L^2(0,L)\to \R$, the
Cauchy problem associated with
\eqref{csystem1} and $f(t)=F(v(t,\cd))$ is
locally well-posed in the sense of the
following theorem, which is proved in the
appendix of this paper.

\begin{theorem}\label{well}
Let $F:L^2(0,L)\to \R$ be continuous linear
map and let $T_0\in (0,+\infty)$. Then, for
every $v^0\in L^2(0,L)$,
 there exists
at most one solution $v\in C^0([0,T_0];
L^2(0,L))\cap L^2(0,T_0;H^1_0(0,L))$ of
\eqref{csystem1} with $f(t)=F(v(t,\cd))$
satisfying the initial condition
$v(0,\cdot)=v^0(\cdot)$.

Moreover, there exists an $r_0>0$ such
that, for every $v^0\in L^2(0,L)$ with
\begin{equation}\label{v0small}
|v^0|_{L^2(0,L)}\leq r_0,
\end{equation}
there exists one
 solution $v\in C^0([0,T_0]; L^2(0,L))\cap
L^2(0,T_0;H^1_0(0,L))$ of \eqref{csystem1}
with $f(t)=F(v(t,\cd))$ satisfying the
initial condition $v(0,\cdot)=v^0(\cdot)$.
\end{theorem}
Let $F:L^2(0,L)\to \R$ be continuous linear
map and let $r_0$ be as in
Theorem~\ref{well} for $T_0\=1$.  Let
$v^0\in L^2(0,L)$ be such that
\eqref{v0small} holds. From
Theorem~\ref{well}, we get that there
exists a unique $T\in (0,+\infty]$ such
that
\begin{itemize}
\item [(i)] There exists $v\in C^0([0,T); L^2(0,L))\cap
L^2_{loc}([0,T);H^1(0,L))$ of
\eqref{csystem1} with $f(t)=F(v(t,\cd))$
satisfying the initial condition
$v(0,\cdot)=v^0(\cdot)$,
\item [(ii)] For every $T'\in (T,+\infty) $ there is no $v\in C^0([0,T'); L^2(0,L))\cap
L^2_{loc}([0,T');H^1(0,L))$ of
\eqref{csystem1} with $f(t)=F(v(t,\cd))$
satisfying the initial condition
$v(0,\cdot)=v^0(\cdot)$.
\end{itemize}
From now on, the $v$ as in (i) will be
called the solution of \eqref{csystem1}
with $f(t)=F(v(t,\cd))$ satisfying the
initial condition $v(0,\cdot)=v^0$. (In
fact \eqref{v0small} is not needed: see
Remark~\ref{largeinitialdata} below.)

It follows from \cite[Lemma 3.5]{Rosier1}
that Problem {\bf (S)} has a negative
answer for the linearized system of
\eqref{csystem1} at $0$ (i.e. the control
system \eqref{csystem1-linear} below) for
every $\lambda
>0$ and every $L\in \cN$ with $\cN$ defined
by
\begin{equation}
\label{defcalN} \cN\=\left\{
2\pi\sqrt{\frac{l^2+lj+j^2}{3}};\,l,j\in\dbZ^+
\right\},
\end{equation}
where $\dbZ^+$ denote the set of positive
integers: $\dbZ^+\=\{1,2,3,\ldots\}$.
(Note, however, that it has been proved in
\cite{CCS} that, for $L=2\pi\in \cN$, $0\in
L^2(0,L)$ is locally asymptotically stable
for our KdV equation \eqref{csystem1} with
$f(\cd)=0$, even if this property is not
true for  the linearized system at $0$ of
our KdV equation \eqref{csystem1} with
$f(\cd)=0$.) In this paper, unless
otherwise specified, we always keep the
assumption that $L\notin\cN$.

The stabilization problems for both linear
and nonlinear systems were studied
extensively in the literature. There are
too many related works to list them
comprehensively here.  As a result,  we
restrict ourselves to KdV equations.   For
KdV equations with internal feedback
controls, we refer to
\cite{KRZ,MM1,MVZ,P,RZ0,RZ1,RZ2,RZ3} and
the rich references therein, as well as the
self-contained survey paper
\cite{2014-Cerpa-MCRF}.  For KdV equations
with boundary feedback controls, we refer
to \cite{CC2,CC1,Zhang1} and the numerous
references therein.

In this paper, we focus on the boundary
feedback stabilization problems. This issue
was first studied in \cite{Zhang1}, where
the author proved that the system
\eqref{csystem1} is locally exponentially
stable in $H_0^1(0,L)$ for $L=1$ and for
the feedback control $f(t)=\a u_x(0,t)$
where $\a\in (0,1)$. In that work, the
exponential decay rate depends on $\a$ but
cannot be arbitrarily large.

There are three natural methods to study
Problem {\bf (S)}. The first one is the so
called ``Gramian approach''. It works well
for time reversible linear systems (see
\cite{1997-Komornik-SICON,1974-Slemrod-SICON,2005-Urquiza-SICON}
for example).  This method was  used in
\cite{CC1} to get the rapid stabilization
for the following linearized system of
\eqref{csystem1}
\begin{equation}\label{csystem1-linear}
\left\{
\begin{array}{ll}\ds
\tilde v_t + \tilde v_x + \tilde v_{xxx} =
0 &\mbox{ in }
(0,+\infty)\times (0,L),\\
\ns\ds \tilde v(t,0) = \tilde v(t,L)=0
&\mbox{ on }
(0,+\infty),\\
\ns\ds \tilde v_x(t,L)=\tilde f(t)&\mbox{
on } (0,+\infty).
\end{array}
\right.
\end{equation}
For every given  $\lambda$, by solving a
suitable linear quadratic regulator
problem, a feedback control law $\tilde
f(t)=\tilde F(\tilde v(t,\cd))$ was
constructed in \cite{CC1} which makes
\eqref{csystem1-linear}  exponentially
stable with an exponential decay rate at
least equal to $\l$. However,  one does not
know how to apply this result to treat the
nonlinear system \eqref{csystem1}.

The second one is the control Lyapunov
function method. It is well known that
Lyapunov functions  are very useful for the
study the asymptotic stability of dynamical
systems. For a control system, with the aid
of a suitable choice of feedback laws,
there are more ``chances'' for a given
function to be a Lyapunov function. Thus,
Lyapunov functions are even more useful for
the stabilization of control systems than
for the stability of dynamical systems
without control (see \cite[Section
12.1]{2007-Coron-book} for example). This
method can be used to get rapid
stabilization for some partial differential
equations (see, for example,
\cite{2008-Coron-Bastin-dAndrea-Nocel-SICON}).
However this method is difficult to use
here, as well as for many other control
systems, to get rapid stabilization. The
reason is that the ``natural'' control
Lyapunov  functions do not lead to
arbitrarily large exponential decay rate.
See, for example, \cite[Remark 12.9, p.
318]{2007-Coron-book}.

The third one is the backstepping method.
It is now a standard method for finite
dimensional control systems (see, e.g.,
\cite{1995-Krstic-Kanellakopoulos-Kokotovic-book},
\cite[pages 242-246]{1998-Sontag-book} and
\cite[Section 12.5]{2007-Coron-book}). The
first adaptations of this method to control
systems modeled by partial differential
equations were given in
\cite{1998-Coron-Andrea-IEEE} and
\cite{2000-Liu-Krstic-SCL}. This method has
been used on the discretization of partial
differential equations in
\cite{2003-Bokovic-Balogh-Krstic-MCSS}.  A
key modification of the method by using a
Volterra transformation is introduced in
\cite{BKL}. This last article has been the
starting point of many works. See, in
particular, the references given in
\cite{KS1}, where  a systematic and clear
introduction to this method is given. Note
that this method can be useful to handle
nonlinearities as it was shown in
\cite{CC2,2013-Coron-Vazquez-Krstic-Bastin-SICON}.
In \cite{CC2}, the authors studied the
rapid stabilization problem with a left
Dirichlet boundary control by the
backstepping method (in that case the
assumption $L\not \in \mathcal{N}$ is no
longer required). They consider the
following system
\begin{equation}\label{csystem2}
\left\{
\begin{array}{ll}\ds
\check v_t + \check v_x + \check v_{xxx} +
\check v\check v_x = 0 &\mbox{ in }
(0,+\infty)\times (0,L),\\
\ns\ds \check v(t,0) = \check f(t), \check
v(t,L)=0 &\mbox{ on }
(0,+\infty), \\
\ns\ds \check v_x(t,L)=0 &\mbox{ on }
(0,+\infty),
\end{array}
\right.
\end{equation}
and proved that, for every $\l>0$, there
exists a continuous linear feedback control
$\check f(t) = \int_0^L \check k(0,y)\check
v(t,y)dy$, an $r>0$ and a $C>0$ such that
$$
|\check v(t,\cd)|_{L^2(0,L)}\leq Ce^{-\l
t}|\check v(0,\cd)|_{L^2(0,L)}
$$
for every solution of \eqref{csystem2}
satisfying $|\check v(0,\cd)|_{L^2(0,L)}
\leq r$. Here the function $$\check k:
\{(x,y):\,x\in [0,L],\,y\in [x,L]\} \to
\mathbb{R}$$ is the solution to
\begin{equation}\label{system2.1}
\left\{
\begin{array}{ll}\ds
\check k_{xxx} + \check k_x + \check k_{yyy} + \check k_y + \l\check k = 0 &\mbox{ in } \{(x,y):\,x\in [0,L],\,y\in [x,L]\},\\
\ns\ds \check k(x,L)=0 &\mbox{ in } [0,L],\\
\ns\ds \check k(x,x)=0 &\mbox{ in } [0,L],\\
\ns\ds \check k_x(x,x)=\frac{\l}{3}(L-x)
&\mbox{ in } [0,L],
\end{array}
\right.
\end{equation}
as in the backstepping approach. %

The main difference between the system
\eqref{csystem1} and the system
\eqref{csystem2} is the position where the
control is acting. It is well known that
the control properties of KdV equations
depend in a crucial way on the location of
the controls. For instance, the system is
only null controllable  if the control acts
on the left Dirichlet boundary condition
and homogeneous data is considered at the
right end point of the interval (see
\cite{Rosier2}). On the other hand, if the
system is controlled from the right
boundary condition, then it is exactly
controllable (see \cite{Rosier1}). In
\cite{CC2}, the authors also pointed out
the difficulties for employing the
backstepping method to solve our rapid
stabilization problem for \eqref{csystem1}.

In this paper, in order to stabilize the
solutions of \eqref{csystem1}, we use a
more general integral transform  on the
state $v$ than the one allowed by  the
backstepping method. In fact, no
restriction is put on the integral
transform which is a priori considered: it
takes the following form
\begin{gather}
\label{transformationv-w}
w(x)\=v(x)-\int_0^L
k(x,y)v(y)dy\=v(x)-(Kv)(x).
\end{gather}
Note that no assumption is made on the
support of $k(\cd,\cd)$. This is in
contrast with the backstepping approach
where the support of $k$ is assumed to be
included in one of the following triangles
$$\{(x,y)\in [0,L]^2; \, x\in [0,L],\, y\in
[0,x]\}, \q \{(x,y)\in [0,L]^2; \, y\in
[0,L],\, x\in [0,y]\}.$$

Now let us briefly explain the idea  for
the linearized KdV equation
\eqref{csystem1-linear}.
One can check that (the computations are
similar but simpler than the ones for the
KdV equation \eqref{csystem1} given in
Section~\ref{secproofth}, we omit them), if
\begin{equation}\label{system3}
\left\{
\begin{array}{ll}\ds
k_{yyy} + k_y  + k_{xxx} + k_x  + \l k = \l
\d(x-y) &\mbox{ in }
(0,L)\times (0,L),\\
\ns\ds k(x,0) = k(x,L) =0 &\mbox{ on } (0,L),\\
\ns\ds k_y(x,0) = k_y(x,L) =0 &\mbox{ on } (0,L),\\
\ns\ds k(0,y) = k(L,y)=0 &\mbox{ on }
(0,L),
\end{array}
\right.
\end{equation}
where $\d(x-y)$ denotes the Dirac measure
on the diagonal of the square $[0,L]\times
[0,L]$, has a solution $k$ which is smooth
enough, then, for every solution $\tilde v$
of \eqref{csystem1-linear}, $\tilde
w(t,x)\=\tilde v(t,x) - \int_0^L
k(x,y)\tilde v(t,y)dy$ is a solution of
\begin{equation}\label{csystem1-linear-w}
\left\{
\begin{array}{ll}\ds
\tilde w_t + \tilde w_x + \tilde w_{xxx}
+\lambda \tilde w = 0 &\mbox{ in }
(0,+\infty)\times (0,L),\\
\ns\ds \tilde w(t,0) = \tilde w(t,L)=0
&\mbox{ on }
(0,+\infty),\\
\ns\ds \tilde w_x(t,L)=\tilde f(t)-\int_0^L
k_x(L,y)\tilde v(t,y)dy&\mbox{ on }
(0,+\infty).
\end{array}
\right.
\end{equation}
If we define the feedback law $F(\cd)$ by
\begin{equation}\label{deffeedback}
\tilde F(\tilde v)=\tilde f(t)\=\int_0^L
k_x(L,y)\tilde v(t,y)dy,
\end{equation}
then the last equation of
\eqref{csystem1-linear-w}  becomes $\tilde
w_x(t,L)=0$. Multiplying the first equation
of \eqref{csystem1-linear-w} by $\tilde w$
and integrating on $[0,L]$ we get, using
integration by parts together with the
boundary conditions of
\eqref{csystem1-linear-w} and $\tilde
w_x(t,L)=0$,
\begin{equation}\label{eqtempw}
\frac{d}{dt}\int_0^L\tilde w(t,x)^2dx \leq
-\lambda \int_0^L\tilde w(t,x)^2dx,
\end{equation}
which shows the exponential decay of the
$L^2$-norm of $\tilde w$. In order to get
the same
 exponential decay of the $L^2$-norm of
$\tilde v$, it suffices to prove that

1. $|v|_{L^2(0,L)}\leq C|v-Kv|_{L^2(0,L)}$
 for some $C>0$  independent of
 $v\in L^2(0,L)$.

\noindent Furthermore, in order to show
that \eqref{est-sur-v} holds, we also need
to prove that

2. $k(\cd,\cd)$ is smooth enough so that
the nonlinearity $vv_x$ will not be a
problem provided that
$|v(0,\cdot)|_{L^2(0,L)}$ is small enough.

 We will check these
two points together with the existence of
$k(\cd,\cd)$ satisfying \eqref{system3} and
therefore we prove the following local
rapid stabilization result.
\begin{theorem}\label{th1}
For every $\l>0$, there exist a continuous
linear feedback control law
$F:L^2(0,L)\to \R$, an $r\in
(0,+\infty)$ and a $C>0$ such that, for
every $v^0\in L^2(0,L)$ satisfying
$|v^0|_{L^2(0,L)}\leq r$, the solution $v$
of \eqref{csystem1} with $f(t)\=
F(v(t,\cdot))$ satisfying the initial condition $v(0,\cdot)=v^0(\cdot)$
 is defined on $[0,+\infty)$
and satisfies
\begin{equation}\label{th1-eq1}
|v(t,\cd)|_{L^2(0,L)}\leq Ce^{-\frac{\l}{2}
t}|v(0,\cd)|_{L^2(0,L)}, \q\mbox{ for every
}\,t \geq 0.
\end{equation}
\end{theorem}

To show the existence and uniqueness of
$k(\cd,\cd)$, we have to utilize the fact
that the linear KdV equation is exact
controllable by Neumann boundary control on
the right end point of the interval. This
is different from the backstepping method,
which does not use any controllability
result.

Although the aim of this paper is to study
the rapid stabilization of KdV equations,
the method introduced here can be applied
to many other partial differential
equations on one dimensional bounded
domains, such as heat equations,
Schr\"{o}dinger equations, beam equations
and Kuramoto-Sivashinsky equations. In
particular,  we have applied it in
\cite{Coron-Lu-KS} to the following
Kuramoto-Sivashinsky control system
\begin{equation}\label{csystemKS}
\left\{
\begin{array}{ll}\ds
v_t + v_{xxxx} + \l v_{xx} + vv_x = 0
&\mbox{
in } (0,+\infty)\times (0,1),\\
\ns\ds v(t,0)=v(t,1)=0 &\mbox{ on } (0,+\infty),\\
\ns\ds v_{xx}(t,0)= f(t),\,v_{xx}(t,1)=0
&\mbox{ on } (0,+\infty),
\end{array}
\right.
\end{equation}
where, at time $t\in [0,+\infty)$, the
state is $v(t,\cdot) \in L^2(0,1)$ and the
control is $f(t)\in \R$.
 In those
cases, the main difference is that  the
equations satisfied by the kernel $k$ are
different from \eqref{system3}. However,
one can follow the same strategy  to show
the existence and regularity of the kernel
and the invertibility of the transform.

The rest of this paper is organized as
follows
\begin{itemize}
\item in Section \ref{secconstructionk}, we
establish the well-posedness of
\eqref{system3} and study the regularity of
its solution,
\item in Section \ref{secinvert}, we prove the invertibility of $I-K$,
\item in Section \ref{secproofth}, we conclude the
proof of Theorem \ref{th1} by using the
results established in
Section~\ref{secconstructionk} and
Section~\ref{secinvert}.
\end{itemize}
%

\section{Well-posedness of \eqref{system3}}
\label{secconstructionk}


This section is devoted to the study of the
equation \eqref{system3}. We first
introduce the definition of the
(transposition) solution to
\eqref{system3}. Let
\begin{equation}\label{def cE}
\begin{array}{ll}\ds
\cE\=\big\{\rho \in C^\infty([0,L]\times [0,L]):\, \rho(0,y)=\rho(L,y)=\rho(x,0)=\rho(x,L)=0,\\
\ns\ds \qq \;
\rho_x(0,y)=\rho_x(0,L)=0\big\}
\end{array}
\end{equation}
and let $\cG$ be the set of $k\in
H_0^1((0,L)\times (0,L))$ such that
\begin{gather}
\label{regkx}
  \left( x\in (0,L)\mapsto k_x(x,\cdot)\in L^2(0,L)\right) \in C^0([0,L];L^2(0,L)),
\\
\label{regky}
  \left( y\in (0,L)\mapsto k_y(\cdot,y)\in L^2(0,L)\right) \in C^0([0,L];L^2(0,L)),
\\
\label{ky0} k_y(\cdot,0)=k_y(\cdot,L)=0
\text{ in } L^2(0,L).
\end{gather}
 We call $k(\cd,\cd)\in \cG$ a solution to \eqref{system3} if
\begin{equation}
\label{deftransposition}
\begin{array}{ll}\ds
\int_0^L\int_0^L \big[\rho_{yyy}(x,y)+\rho_y(x,y) + \rho_{xxx}(x,y) + \rho_x(x,y) - \l \rho(x,y) \big]k(x,y)dxdy\\
\ns\ds + \int_0^L \l \rho(x,x)dx=0,
\q\mbox{ for every } \rho\in\cE.
\end{array}
\end{equation}

We have the following well-posedness result
for \eqref{system3}.
\begin{lemma}\label{well-posed lm}
Let  $\l\neq 0$. Equation \eqref{system3}
has one and only one solution in $\cG$.
\end{lemma}
{\it Proof of Lemma~\ref{well-posed lm}}\,:
The proof is divided into two steps:
\begin{itemize}
\item Step 1: proof of the uniqueness of the solution to
\eqref{system3};
\item Step 2: proof of the existence of a solution to \eqref{system3} with the required regularity.
\end{itemize}

\vspace{0.3cm}

\begin{center}
\textbf{Step 1: proof of the uniqueness of
the solution to \eqref{system3}}
\end{center}

Assume that $k_1(\cd,\cd)$ and
$k_2(\cd,\cd)$ are two solutions of
\eqref{system3}. Let
$k_3(\cd,\cd)\=k_1(\cd,\cd)-k_2(\cd,\cd)$.
Then $k_3(\cd,\cd)$ is a (transposition:
\eqref{deftransposition} holds without the
last integral term) solution of
\begin{equation}\label{system4}
\left\{
\begin{array}{ll}\ds
k_{3,yyy} + k_{3,y}  + k_{3,xxx} + k_{3,x}
+ \l k_3 = 0 &\mbox{ in }
(0,L)\times (0,L),\\
\ns\ds k_3(x,0) = k_3(x,L) =0 &\mbox{ on } (0,L),\\
\ns\ds k_{3,y}(x,0) = k_{3,y}(x,L) =0 &\mbox{ on } (0,L),\\
\ns\ds k_3(0,y) = k_3(L,y)=0 &\mbox{ on }
(0,L).
\end{array}
\right.
\end{equation}

Let us define an unbounded linear operator
$A: D(A)\subset L^2(0,L)\to L^2(0,L)$ as
follows.
\begin{equation}\label{DA}
\left\{
\begin{array}{ll}\ds
D(A)\=\{ \f:\,\f\in
H^3(0,L),\,\f(0)=\f(L)=0,\,\f_x(0)=\f_x(L)
\},\\
\ns\ds A\f \=-\f_{xxx}-\f_x.
\end{array}
\right.
\end{equation}
The operator $A$ is a skew-adjoint operator
with compact resolvent. Furthermore, since
$L\notin\cN$, we have $L\not \in 2\pi
\mathbb{Z}^+$, which, as one easily checks,
implies that $0$ is not an eigenvalue of
$A$. Denote by
$\{i\mu_j\}_{j\in\dbZ\setminus\{0\} }$,
$\mu_j\in \R$, the eigenvalues of $A$,
which are organized in the following way:
\begin{gather}
\label{ordreeig} \cds\leq  \mu_{-2}\leq
\mu_{-1} < 0 <  \mu_1 \leq \mu_2 \leq \cds
.
\end{gather}
Let us point out that all these eigenvalues
are simple. Indeed, assume that $\mu_j$ is
of multiplicity at least $2$ and $\f_{1,j}$
and $\f_{2,j}$ are two linear independent
eigenfunctions corresponding to $i\mu_j$.
Let
\begin{equation}
\f_j \= \left\{
\begin{array}{ll}
\ds
\f_{2,j}-\frac{\f_{2,j}'(0)}{\f_{1,j}'(0)}\f_{1,j}&\text{
if } \f_{1,j}'(0)\not = 0,
\\
\f_{1,j} &\text{ if } \f_{1,j}'(0)= 0.
\end{array}
\right.
\end{equation}
Since $\f_{1,j}$ and $\f_{2,j}$ are
linearly independent, we conclude that
$\f_j\neq 0$. Then $\f_j$  is  a nonzero
eigenfunction  for this eigenvalue $i\mu_j$
such that
\begin{equation}
\label{tout0}
\f_j(0)=\f_j(L)=\f_j'(0)=\f_j'(L)=0.
\end{equation}
By \cite[Lemma 3.5]{Rosier1}, our
assumption $L\notin \mathcal{N}$ implies
that,  for every a nonzero eigenfunction
$\f_j$  of $A$,
\begin{equation}\label{fjx(0)not0}
\f_{j}'(0)\not=0,
\end{equation}
which is in contradiction with
\eqref{tout0}. Hence all these eigenvalues
are simple and all the inequalities in
\eqref{ordreeig} are strict inequalities.
Let us write
$\{\f_j\}_{j\in\dbZ\setminus\{0\} }$ for
the corresponding eigenfunctions with
$|\f_j|_{L^2(0,L)}=1$
($j\in\dbZ\setminus\{0\} $). It is well
known that
$\{\f_j\}_{j\in\dbZ\setminus\{0\} }$
constitutes an orthonormal basis of
$L^2_{\mathbb{C}}(0,L)$.

In what follows, for simplicity, in
$L^2_{\mathbb{C}}(0,L)$,
$L^\infty_{\mathbb{C}}(0,L)$ etc., when no
confusion is possible,    we omit the index
$\mathbb{C}$.

Next, we analyze the asymptotic behavior of
$\mu_j$ as $|j|\rightarrow +\infty$. This
is done in \cite{CC1}. However,  we present
it also here since we need similar
information  for other sequences of numbers
which will be introduced later on. To this
end, we first consider the following
boundary problem,  with $\mu$ given in
$\mathbb{R}$,
\begin{equation}\label{system17}
\left\{
\begin{array}{ll}\ds
\f''' + \f' + i\mu\f=0 &\mbox{ in
} [0,L],\\
\ns\ds \f(0)=\f(L)=0,\\
\ns\ds \f'(0)=\f'(L).
\end{array}
\right.
\end{equation}
To solve the equation \eqref{system17}, we
study the following algebraic equation in
$\mathbb{C}$
\begin{equation}\label{system18}
s^3+s+i\mu=0.
\end{equation}
Denote by $s_1$, $s_2$ and $s_3$ the three
roots of \eqref{system18}. These roots  are
distinct if and only if $4\neq 27\mu^2$,
which we assume from now on. For the
special case $4= 27\mu^2$, the arguments
given below can be adapted and we omit this
adaptation.

Since the function $x\in \mathbb{R}\mapsto
2x(4x^2-1)\in \mathbb{R}$ is a surjective
map from $\dbR$ to $\dbR$, we get that
there always exists at least one
$\tau\in\dbR$ such that
\begin{equation}
\label{detau} \mu = 2\tau(4\tau^2-1).
\end{equation}
 By this,  the
three roots of \eqref{system18} are
$$
s_1=\sqrt{ 3\tau^2-1 }-i\tau,\;s_2=-\sqrt{
3\tau^2-1 }-i\tau,\;s_3 = 2i\tau.
$$
Since we want to analyze the asymptotic
behavior of $\mu_j$ as $|j|\rightarrow
+\infty$, we may assume that $|\tau|$ is
large enough.  Thus, we only consider the
case $3\tau^2>1$. Now the eigenfunction  is
$$
\f(x) = e^{-i\tau x}\big[
r_1\cosh(\sqrt{3\tau^2-1}x) +
r_2\sinh(\sqrt{3\tau^2-1}x) \big] +
r_3e^{2i\tau x}
$$
for some suitable complex numbers $r_1$,
$r_2$ and $r_3$ which are not all equal to
$0$ . The boundary  conditions of
\eqref{system17} are fulfilled if and only
if these three complex numbers satisfy that
$$
\left\{
\begin{array}{ll}\ds
r_1 + r_3 = 0,\\
\ns\ds e^{-i\tau L}\big[
r_1\cosh(\sqrt{3\tau^2-1}L) +
r_2\sinh(\sqrt{3\tau^2-1}L) \big] +
r_3e^{2i\tau L}=0,\\
\ns\ds \q-i\tau r_1 +
2i\tau r_3+r_2\sqrt{3\tau^2-1}\\
\ns\ds = -i\tau e^{-i\tau L} \big[
r_1\cosh(\sqrt{3\tau^2-1}L) +
r_2\sinh(\sqrt{3\tau^2-1}L) \big] +
2i\tau r_3 e^{2i\tau L}\\
\ns\ds \q + e^{-i\tau
L}\sqrt{3\tau^2-1}\big[
r_1\sinh(\sqrt{3\tau^2-1}L) +
r_2\cosh(\sqrt{3\tau^2-1}L) \big].
\end{array}
\right.
$$
Thus, we have that
\begin{equation}\label{9.27-eq1}
\sqrt{3\tau^2\!-\!1}\cos(2\tau
L)\!-\!3\tau\sin(\tau L)
\sinh(\sqrt{3\tau^2\!-\!1}L)\!-
\sqrt{3\tau^2\!-\!1}\cos(\tau
L)\cosh(\sqrt{3\tau^2\!-\!1}L) =0.
\end{equation}
As $|\tau|\rightarrow +\infty$,
$$
\sinh(\sqrt{3\tau^2-1}L)\sim
\cosh(\sqrt{3\tau^2-1}L)\sim \frac{1}{2}
e^{\sqrt{3\tau^2-1}L}\sim \frac{1}{2}
e^{\sqrt{3}|\tau |L}.
$$
Hence, we find that
\begin{equation}\label{system19}
\begin{array}{ll}\ds
e^{\sqrt{3\tau^2-1}L}= \frac{\cos (2\tau
L)}{2\cos(\tau L - \frac{\pi}{3})} + O(1)\,
\text{ as } \tau \to +\infty,\\
\ns\ds e^{\sqrt{3\tau^2-1}L}= \frac{\cos
(2\tau L)}{2\cos(\tau L + \frac{\pi}{3})} +
O(1)\, \text{ as } \tau \to -\infty.
\end{array}
\end{equation}
For $j\in \dbZ $ with $|j|$ large enough,
one checks that there exists a unique
solution $\tau=\tau_j$ of \eqref{9.27-eq1}
in each interval
$[\frac{j\pi}{L},\frac{(j+1)\pi}{L})$, and
one has that
\begin{gather}
\label{estimatetauj} \left\{
\begin{array}{ll}\ds
\tau_j = \frac{j\pi}{L} + \frac{5\pi}{6L} +
O\(\frac{1}{j}\)\,&\mbox{ as
}\;  j  \to + \infty,\\
\ns\ds \tau_j = \frac{j\pi}{L} +
\frac{\pi}{6L} + O\(\frac{1}{j}\)\,&\mbox{
as }\;  j  \to - \infty.
\end{array}
\right.
\end{gather}
 From this and \eqref{detau} we can see that
\begin{equation}\label{5.22-eq1}
\left\{
\begin{array}{ll}\ds
\mu_j = \(\frac{2j\pi}{L}\)^3 +
\frac{40\pi^3}{3}j^2 + O(j)& \mbox{
as }\; j \to + \infty,\\
\ns\ds \mu_j = \(\frac{2j\pi}{L}\)^3 +
\frac{8\pi^3}{3}j^2+ O(j)& \mbox{ as }\; j
\to -\infty.
\end{array}
\right.
\end{equation}
The corresponding eigenfunction $\f_j$ with
$|\f_j|_{L^2(0,L)}=1$ reads
\begin{equation}\label{5.22-eq2}
\f_j(x)\!= \! \a_j \bigg\{\!e^{-i\tau_j
x}\bigg[\!
\cosh\big(\sqrt{3\tau_j^2\!-\!1}x\big)
\!+\! \frac{e^{3i\tau_j
L}-\cosh\big(\sqrt{3\tau_j^2\!-\!1}L\big)}{\sinh\big(\sqrt{3\tau_j^2-1}L\big)}\sinh\big(\sqrt{3\tau_j^2\!-\!1}x\big)
\!\bigg]\! - \!e^{2i\tau_j x}\!\bigg\},
\end{equation}
where $\a_j\in (0,+\infty)$ is a suitable
chosen positive real number such that
$|\f_j|_{L^2(0,L)}=1$.

Let us analyze the asymptotic behavior of
$\a_j$ as $|j|\to +\infty$.  For this, with
some simple computations, we see that
\begin{equation}\label{5.15-eq0}
\!\!\begin{array}{ll}\ds \f_j(x) \3n& \ds=
\a_j \Bigg\{\frac{
\cosh\!\big(\!\sqrt{3\tau_j^2\!-\!\!1}x\big)
\sinh\!\big(\!\sqrt{3\tau_j^2\!-\!\!1}L\big)\!\!-\!\sinh\!\big(\!\sqrt{3\tau_j^2\!-\!1}x\big)\cosh\!\big(\!\sqrt{3\tau_j^2\!-\!\!1}L\big)}{\sinh\big(\sqrt{3\tau_j^2-1}L\big)}
e^{-i\tau_j x}
\\
\ns&\ds \q  + \frac{e^{3i\tau_j
L}\sinh\big(\sqrt{3\tau_j^2-1}x\big)}{\sinh\big(\sqrt{3\tau_j^2-1}L\big)}e^{-i\tau_j
x} - e^{2i\tau_j x}\Bigg\}.
\end{array}
\end{equation}
Then, we find that, as $|j|\to +\infty$,
\begin{equation}\label{5.15-eq1}
\begin{array}{ll}\ds
\Bigg(\int_0^L\Bigg|\f_j(x)- \a_j \Bigg[
\frac{2e^{3i\tau_j
L}\sinh\big(\sqrt{3\tau_j^2 - 1}x\big) +
e^{\sqrt{3\tau_j^2 -
1}(L-x)}}{2\sinh\big(\sqrt{3\tau_j^2-1}L\big)}e^{-i\tau_j
x}- e^{2i\tau_j x}\Bigg]
\Bigg|^2dx\Bigg)^{\frac{1}{2}}\\
\ns\ds=
O\big(\a_je^{-\sqrt{3\tau_j^2-1}L}\big),
\end{array}
\end{equation}
which, using $|\f_j|_{L^2(0,L)}=1$, gives
us
\begin{equation}\label{6.29-eq1}
\a_j \to \frac{1}{\sqrt{L}} \mbox{ as }
|j|\to + \infty.
\end{equation}

We now estimate $\f_{j,x}$. From
\eqref{5.22-eq2}, we obtain that
\begin{equation}\label{8.30-eq14}
\begin{array}{ll}\ds
\f_{j,x}(x)  = \a_j \Bigg\{-i\tau_j
e^{-i\tau_j x}\Bigg[
\cosh\big(\sqrt{3\tau_j^2\!-\!1}x\big)
\!+\! \frac{e^{3i\tau_j
L}-\cosh\big(\sqrt{3\tau_j^2\!-\!1}L\big)}{\sinh\big(\sqrt{3\tau_j^2-1}L\big)}\sinh\big(\sqrt{3\tau_j^2\!-\!1}x\big)
\Bigg]\!\\
\ns \ds \qq\; +e^{-i\tau_j x} \Bigg[\!
\sqrt{3\tau_j^2\!\!-\!1}\sinh\big(\sqrt{3\tau_j^2\!\!-\!1}x\big)
\!+\! \frac{e^{3i\tau_j
L}\!\!-\!\cosh\big(\sqrt{3\tau_j^2\!\!-\!1}L\big)}{\sinh\big(\sqrt{3\tau_j^2-1}L\big)}\!\sqrt{3\tau_j^2\!-\!1}\cosh\!\big(\!\sqrt{3\tau_j^2\!\!-\!1}x\big)
\Bigg]\\
\ns\ds \qq\; - 2i\tau_je^{2i\tau_j x}
\Bigg\}.
\end{array}
\end{equation}
 From \eqref{6.29-eq1} and \eqref{8.30-eq14}, we get that, for $|j|$
large enough,
\begin{equation}\label{8.30-eq15}
\begin{array}{ll}\ds
\q|\f_{j,x}|_{L^\infty(0,L)}\\
\ns\ds \leq  \max_{x\in [0,L]}\Bigg|\a_j
\tau_j \Bigg[ \frac{2e^{3i\tau_j
L}\sinh\big(\sqrt{3\tau_j^2\!-\!1}x\big)+e^{\sqrt{3\tau_j^2-1}L-\sqrt{3\tau_j^2-1}x}
-
e^{-\sqrt{3\tau_j^2-1}L+\sqrt{3\tau_j^2-1}x}}{2\sinh\big(\sqrt{3\tau_j^2-1}L\big)}
\Bigg]\Bigg| \\
\ns\ds \q \!+ \max_{x\in
[0,L]}\Bigg|\a_j\sqrt{3\tau_j^2\!\!-\!1}\Bigg[
\frac{2e^{3i\tau_j
L}\!\sinh\big(\sqrt{3\tau_j^2\!\!-\!1}x\big)\!\!-\!e^{\sqrt{3\tau_j^2\!-1}L-\sqrt{3\tau_j^2\!-1}x}
\!\!\!-\!
e^{-\sqrt{3\tau_j^2\!-1}L+\sqrt{3\tau_j^2\!-1}x}}{2\sinh\big(\sqrt{3\tau_j^2-1}L\big)}
\Bigg] \Bigg| \\
\ns\ds \q + 2|\tau_j|\\
\ns\ds \leq C|\tau_j|,
\end{array}
\end{equation}
From \eqref{estimatetauj} and
\eqref{8.30-eq15}, one gets, for every
$j\in \dbZ^+$,
\begin{equation}\label{estimateLinftyfjx(x)}
  |\f_{j,x}|_{L^\infty(0,L)}\leq C|j|.
\end{equation}
 From \eqref{8.30-eq14}, as $|j|\to +\infty$,
\begin{equation}\label{8.30-eq5}
\begin{array}{ll}\ds
\f_{j,x}(0)\3n&\ds = \a_j \Bigg[
\frac{e^{3i\tau_j
L}-\cosh\big(\sqrt{3\tau_j^2-1}L\big)}{\sinh\big(\sqrt{3\tau_j^2-1}L\big)}\sqrt{3\tau_j^2-1}
- 3i\tau_j \Bigg]\\
\ns&\ds = \a_j\(-3i\tau_j - \sqrt{3\tau_j^2
- 1} \, \)+ \a_j O\big(e^{-\sqrt{3\tau_j^2
- 1}L} \big)
\end{array}
\end{equation}
and
\begin{equation}\label{8.30-eq6}
\begin{array}{ll}\ds
\f_{j,x}(L)\3n&\ds = \a_j
\Bigg\{-3i\tau_je^{2i\tau_j L}+e^{-i\tau_j
L}\Bigg[\sqrt{3\tau_j^2-1}
\sinh\big(\sqrt{3\tau_j^2-1}L\big)\\
\ns&\ds\q  +
\sqrt{3\tau_j^2-1}\frac{e^{3i\tau_j
L}-\cosh\big(\sqrt{3\tau_j^2-1}L\big)}{\sinh\big(\sqrt{3\tau_j^2-1}L\big)}\cosh\big(\sqrt{3\tau_j^2-1}L\big)
\Bigg]\Bigg\}\\
\ns&\ds = \a_j\(-3i\tau_j + \sqrt{3\tau_j^2
- 1} \, \)e^{2i\tau_j L} + \a_j
O\big(e^{-\sqrt{3\tau_j^2 - 1}L}\,\big).
\end{array}
\end{equation}

Let us write
\begin{gather}
\label{k3=}
k_3(x,y)=\sum_{j\in\dbZ\setminus\{0\}}
\psi_j(x)\f_j(y)
\end{gather}
 for the solution to
\eqref{system4}. Then, we get that $\psi_j$
solves
\begin{equation}\label{system16.1}
\left\{
\begin{array}{ll}\ds
\psi_j''' + \psi_j' + \l \psi_j - i\mu_j
\psi_j = 0
&\mbox{ in }(0,L),\\
\ns\ds \psi_j(0) = \psi_j(L)=0.
\end{array}
\right.
\end{equation}
Let $c_j \= \psi_{j,x}(L) - \psi_{j,x}(0)$
($j\in \mathbb{Z}\setminus\{0\}$). We
consider the following equation:
\begin{equation}\label{system16}
\left\{
\begin{array}{ll}\ds
\check\psi_j''' + \check\psi_j' + \l
\check\psi_j - i\mu_j \check\psi_j = 0
&\mbox{ in }(0,L),\\
\ns\ds \check\psi_j(0) = \check\psi_j(L)=0,\\
\ns\ds \check\psi_j'(L) -
\check\psi_j'(0)=1.
\end{array}
\right.
\end{equation}
Since $\l-i\mu_j$ is not an eigenvalue of
$A$ (recall that $\l\neq 0$), the equation
\eqref{system16} has one and only one
solution. Moreover
\begin{equation}\label{10.5-eq9}
\psi_j=c_j\check\psi_j \mbox{ for every
}j\in\dbZ\setminus\{0\}.
\end{equation}
Denote by $r^{(1)}_j$, $r^{(2)}_j$ and
$r^{(3)}_j$ the roots of
\begin{equation}\label{3.27-eq1}
r^3 + r=-\l + i\mu_j.
\end{equation}
Let $\si_j\in\dbR$ be such that
\begin{equation}\label{sigmaj}
2\si_j(4\si_j^2-1)= \left\{
\begin{array}{ll}\ds
-\sqrt{\l^2+\mu_j^2}, &\mbox{ if } \mu_j>0,\\
\ns\ds \sqrt{\l^2+\mu_j^2}, &\mbox{ if }
\mu_j< 0.
\end{array}
\right.
\end{equation}
 Then, there exists $C>0$ such that
\begin{equation}
\label{6.17-eq2} |\si_{-j} - \tau_j|\leq
\frac{C}{1+j^2} \;\mbox{ and } \;
|\si_{-j}^2 - \tau_j^2|\leq \frac{C}{1+|j|}
, \, \forall\, j\in\dbZ\setminus\{0\}.
\end{equation}
Let $\hat\si_j\in\dbC$  be such that
\begin{equation}\label{6.20-eq1}
2\hat\si_j(4\hat\si_j^2-1)=-i\l -\mu_j.
\end{equation}
Roughly speaking, there are three complex
numbers which satisfy \eqref{6.20-eq1}. For
a good choice of $\hat\si_j$, we have
$$
\lim_{|j|\to +\infty} |\si_j-\hat\si_j|=0.
$$
In this case, it is easy to see that  there
is a constant $C>0$ such that
\begin{equation}\label{6.20-eq2}
|\si_j-\hat \si_j| \leq \frac{C}{1+j^2}, \,
\forall j\in \mathbb{Z}\setminus\{0\}.
\end{equation}
The solutions to \eqref{3.27-eq1} read
$$
r_j^{(1)}=\sqrt{ 3\hat\si_j^2-1
}-i\hat\si_j,\;r_j^{(2)}=-\sqrt{
3\hat\si_j^2-1 }-i\hat\si_j,\;r_j^{(3)} =
2i\hat\si_j.
$$

 For every $j\in \mathbb{Z}\setminus\{0\}$, let us define
\begin{gather}
\label{defkappa} \k_j^{(1)}\=
\frac{-i\hat\si_j + \sqrt{ 3\hat\si_j^2-1 }
- \sqrt{ 3 \si_j^2-1 }}{-i\si_j}, \q
\k_j^{(2)}\= \frac{\hat\si_j}{\si_j},
\\
\label{deftheta} \th_j(x) \=
e^{(-i\hat\si_j - \sqrt{3\hat\si_j^2-1})x}
- e^{(-i\hat\si_j+\sqrt{3\hat \si_j^2-1} -
2\sqrt{3\si_j^2-1})x},\, x\in [0,L].
\end{gather}
Let $J\in \mathbb{Z}^+$ be such that
\begin{gather}
\label{defpropertyJ} \left(j\in \mathbb{Z}
\text{ and } |j|\geq J\right)\Rightarrow
\left(3 \si_j^2-1\geq 1 \text{ and }
|3\hat\si_j^2-1|\geq 1 \right).
\end{gather}
(The existence of such a $J$ follows from
\eqref{estimatetauj}, \eqref{6.17-eq2} and
\eqref{6.20-eq2}). From  \eqref{6.20-eq2}
and \eqref{deftheta}, one gets the
existence of $C_1>0$ and $C>0$ such that,
for every $x\in [0,L]$ and every
$j\in\mathbb{Z}\setminus\{0\}$ with
$|j|\geq J$,
\begin{equation}\label{6.20-eq3}
\begin{array}{ll}\ds
|\th_j(x)|\3n&\ds\leq \Big|e^{(-i\hat\si_j
- \sqrt{3\hat\si_j^2-1}+\sqrt{3
\si_j^2-1})x} -e^{(-i\hat\si_j +
\sqrt{3\hat\si_j^2-1}
-\sqrt{3 \si_j^2-1})x}\Big|e^{-\sqrt{3 \si_j^2-1}x}\\
\ns&\ds \leq \Big|e^{(-i\hat\si_j -
\sqrt{3\hat\si_j^2-1}+\sqrt{3
\si_j^2-1})x}\Big|\Big|1-e^{2(
\sqrt{3\hat\si_j^2-1}
-\sqrt{3 \si_j^2-1})x}\Big|e^{-\sqrt{3 \si_j^2-1}x} \\
\ns&\ds \leq
C_1e^{\frac{C_1}{j^2}}\frac{C_1}{j^2}e^{-\sqrt{3
\si_j^2-1}x}
\leq \frac{C}{j^2}e^{-\sqrt{3 \si_j^2-1}x}.\\
\end{array}
\end{equation}
Similarly, there exists $C>0$ such that,
for every $x\in [0,L]$ and all integer $j$
with $|j|\geq J$,
\begin{gather}
\label{estimatethetaprime} |\th_j'(x)|\leq
\frac{C\sqrt{3 \si_j^2-1}}{j^2}e^{-\sqrt{3
\si_j^2-1}x}.
\end{gather}
From \eqref{6.20-eq3} and
\eqref{estimatethetaprime}, one gets the
existence of $C>0$ such that, for all
integer $j$ with $|j|\geq J$,
\begin{equation}\label{6.20-eq4}
|\th_j(L)|\leq Ce^{-\sqrt{3 \si_j^2-1}L},
\q |\th_j'(L)|\leq C\sqrt{3
\si_j^2-1}e^{-\sqrt{3 \si_j^2-1}L}.
\end{equation}

Similar to the argument for the analysis of
$\mu_j$ and $\f_j$, we can get that, if
$3\si_j^2\leq 1$,
$$
\begin{array}{ll}\ds
\check\psi_j(x)\3n&\ds =
e^{-i\k_j^{(1)}\si_j x}\bigg[
\b_j^{(1)}\cos\big(\sqrt{1 -
3\si_j^2}x\big) +
\b_j^{(2)}\sin\big(\sqrt{1 -
3\si_j^2}x\big)
\bigg] + \b_j^{(3)}e^{2i\k_j^{(2)}\si_j x} \\
\ns&\ds \q  +
\big(\b_j^{(1)}-\b_j^{(2)}\big)\th_j(x),
\end{array}
$$
and, if $3\si_j^2>1$,
$$
\begin{array}{ll}\ds
\check\psi_j(x) \3n&\ds =
e^{-i\k_j^{(1)}\si_j x}\bigg[
\b_j^{(1)}\cosh\big(\sqrt{3\si_j^2 -
1}x\big) +
\b_j^{(2)}\sinh\big(\sqrt{3\si_j^2 -
1}x\big) \bigg] +
\b_j^{(3)}e^{2i\k_j^{(2)}\si_j x}
\\
\ns&\ds \q + \big(\b_j^{(1)} -
\b_j^{(2)}\big)\th_j(x).
\end{array}
$$

From now on, we assume that $|j|$ is large
enough so that $3\si_j^2>1$. From the
boundary conditions satisfied by
$\check\psi_j$ -see \eqref{system16}-, we
get
\begin{equation}\label{3.29-eq2}
\!\left\{
\begin{array}{ll}\ds
\!\b_j^{(1)} + \b_j^{(3)} = 0,\\
\ns\ds \!e^{-i\k_j^{(1)}\si_j L}\bigg[
\b_j^{(1)}\cosh\big(\sqrt{3\si_j^2-1}L\big)
+
\b_j^{(2)}\sinh\big(\sqrt{3\si_j^2-1}L\big)
\bigg] +
\b_j^{(3)}e^{2i\k_j^{(2)}\si_j L}\\
\ns\ds   \!+ \big(\b_j^{(1)}\!-\!\b_j^{(2)}\big)\th_j(L)=0,\\
\ns\ds \q\!-i\k_j\si_j\b_j^{(1)} +
2i\k_j\si_j\b_j^{(3)}+\b_j^{(2)}\sqrt{3\si_j^2-1}
+ 1\\
\ns\ds \!= -i\k_j\si_j e^{-i\k_j^{(1)}\si_j
L} \!\bigg[
\b_j^{(1)}\cosh\big(\sqrt{3\si_j^2-1}L\big)
\!+\!
\b_j^{(2)}\sinh\big(\sqrt{3\si_j^2-1}L\big)
\bigg]\!\! +\!
2i\si_j\b_j^{(3)} e^{2i\k_j^{(2)}\si_j L}\\
\ns\ds \q\! + e^{-i\k_j^{(1)}\si_j
L}\!\sqrt{3\si_j^2\!-\!1}\bigg[
\b_j^{(1)}\sinh\!\big(\!\sqrt{3\si_j^2\!-\!\!1}L\big)
\!+\!
\b_j^{(2)}\cosh\!\big(\!\sqrt{3\si_j^2\!-\!\!1}L\big)
\bigg]\!+\!\big(\b_j^{(1)}\!\!-\!\b_j^{(2)}\big)\th_j'(L).
\end{array}
\right.
\end{equation}
Then, we know
\begin{equation}\label{5.23-eq1}
\begin{array}{ll}\ds
\check\psi_j(x)\3n &\ds\!= \!\b_j
\Bigg\{\!e^{-i\k_j^{(1)}\si_j x}\Bigg[\!
\cosh\big(\sqrt{3\si_j^2\!-\!1}x\big)\! +\!
\frac{e^{i(\k_j^{(1)}\!+ 2\k_j^{(2)})\si_j
L}\!\!-\!\cosh\big(\sqrt{3\si_j^2\!-\!1}L\big)\!-\!e^{i\k_j^{(1)}\si_j
L}\th_j(L)}{\sinh\big(\sqrt{3\si_j^2-1}L\big)- e^{i\k_j^{(1)}\si_j L}\th_j(L)}\\
\ns&\ds \q\times
\bigg(\sinh\big(\sqrt{3\si_j^2-1}x\big)
-e^{i\k_j\si_j x}\th_j(x) \bigg) \Bigg] -
e^{2i\k_j^{(2)}\si_j x}\Bigg\},
\end{array}
\end{equation}
where $\b_j$ is a suitable chosen complex
number such that the third equality in
\eqref{3.29-eq2} holds.

From \eqref{5.23-eq1}, we have, as $|j|\to
+\infty$,
\begin{equation}\label{10.5-eq2}
\!\!\!\!\!\!\begin{array}{ll}\ds
\check\psi_{j,x}(0)\3n&\ds = \!\b_j \Bigg[-
i \k_j^{(1)}\si_j \!-\! 2i \k_j^{(2)}\si_j
\!+\! \frac{e^{i( \k_j^{(1)}+2
\k_j^{(2)})\si_j
L}\!-\!\cosh\big(\sqrt{3\si_j^2\!-\!1}L\big)\!-\!e^{i\k_j^{(1)}\si_j
L}\th_j(L)}{\sinh\big(\sqrt{3\si_j^2-1}L\big)- e^{i \k_j^{(1)}\si_j L} \th_j(L)}\\
\ns&\ds\q\times\bigg(2\sqrt{3\hat\si_j^2-1}- \sqrt{3 \si_j^2-1}\bigg)\Bigg]\\
\ns&\ds = \b_j\bigg(\!- i \k_j^{(1)}\si_j -
2i \k_j^{(2)}\si_j -
2\sqrt{3\hat\si_j^2-1}+ \sqrt{3 \si_j^2-1}
\,\bigg) + \g_j O\bigg(e^{-\sqrt{3\si_j^2 -
1}L}\bigg)
\end{array}
\end{equation}
and
\begin{equation}\label{10.5-eq3}
\!\!\!\begin{array}{ll}\ds
\check\psi_{j,x}(L)\3n&\ds = \!\b_j
\Bigg\{\!\!-\!i \k_j^{(1)}\si_je^{2i
\k_j^{(2)}\si_j L}\! -\!2i
\k_j^{(2)}\si_je^{2i \k_j^{(2)}\si_j
L}\!+\!e^{-i \k_j^{(1)}\si_j
L}\!\Bigg[\!\sqrt{3\si_j^2\!- \!1}
\sinh\big(\sqrt{3\si_j^2\!- \!1}L\big)\\
\ns&\ds\q \! + \frac{e^{i( \k_j^{(1)}+2
\k_j^{(2)})\si_j
L}\!\!-\!\cosh\big(\sqrt{3\si_j^2\!-\!1}L\big)\!-\!e^{i\k_j^{(1)}\si_j
L}\th_j(L)}{\sinh\big(\sqrt{3\si_j^2-1}L\big)+
e^{i \k_j^{(1)}\si_j L} \th_j(L)}\(\! \sqrt{3\si_j^2\!\!-\!1}\cosh\!\big(\!\sqrt{3\si_j^2\!\!-\!1}\!L \big)\\
\ns&\ds \q -i\k_j^{(1)}\si_j e^{i
\k_j^{(1)}\si_j L} \th_j(L)-i e^{i
\k_j^{(1)}\si_j L}\th_j'(L)\)
\Bigg]\Bigg\}\\
\ns&\ds = \!\b_j\(- i \k_j^{(1)}\si_j - 2i
\k_j^{(2)}\si_j + \sqrt{3\si_j^2 -
1}\,\)e^{2 i \k_j^{(2)}\si_j L}  + \b_j
O\big(e^{-\sqrt{3\si_j^2 - 1}L}\big).
\end{array}
\end{equation}

Since $\f_{j,y}(0)=\f_{j,y}(L)$, from
\eqref{8.30-eq5} and \eqref{8.30-eq6}, we
get
\begin{equation}\label{10.5-eq4}
-3i\tau_j - \sqrt{3\tau_j^2 - 1}+
O\big(e^{-\sqrt{3\tau_j^2 - 1}L}\big)=
\(-3i\tau_j + \sqrt{3\tau_j^2 - 1}\,
\)e^{2i\tau_j L}.
\end{equation}
This, together with the choice of $\si_j$,
$ \k_j^{(1)}$ and $\k_j^{(2)}$, implies
that
\begin{equation}\label{10.5-eq5}
\begin{array}{ll}\ds
\q - i \k_j^{(1)}\si_j - 2i \k_j^{(2)}\si_j - 2\sqrt{3\si_j^2-1}+\sqrt{3\hat\si_j^2-1}\\
\ns\ds = \(- i \k_j^{(1)}\si_j - 2i
\k_j^{(2)}\si_j + \sqrt{3\si_j^2 -
1}\,\)e^{2 i \k_j^{(2)}\si_j L}  +
O\(\frac{1}{j^2}\), \q\mbox{ as } |j|\to
+\infty.
\end{array}
\end{equation}
Thus, we see
\begin{equation}\label{10.5-eq6}
\b_j = O(j^2) \;\mbox{ as } |j|\to +\infty.
\end{equation}

 Let, for $j\in\dbZ\setminus\{0\}$,
\begin{equation}
\label{10.5-eq1} \hat\psi_{j}(x)\= \left\{
\begin{array}{lll}
\ds
\frac{\a_{-j}}{\b_j}\check\psi_{j}(x),&\q
x\in [0,L] &\text{ if }3\sigma_j^2 >1,
\\
\ds \check\psi_j(x) ,&\q x\in [0,L] &\text{
if }3\sigma_j^2 \leq 1.
\end{array}
\right.
\end{equation}
Now, we are going to prove that
$\{\hat\psi_j\}_{j\in\dbZ\setminus\{0\} }$
is a Riesz basis of $L^2(0,L)$.

First, we analyze the asymptotic behavior
of $\hat\psi_j(\cd)$ as $|j|\to +\infty$.
It is clear that
\begin{equation}\label{5.15-eq2}
\begin{array}{ll}\ds
\!\hat\psi_j(x) \3n\!& \ds= \!\a_{-j}
\Bigg\{
\frac{e^{i(\k_j^{(1)}+2\k_j^{(2)})\si_j
L}\sinh\big(\sqrt{3\si_j^2\!-\!1}x\big)}{\sinh\big(\sqrt{3\si_j^2\!-\!1}L\big)\!\!-
\! e^{i\k_j^{(1)}\si_j
L}\th_j(L)}e^{-i\k_j^{(1)}\si_j x}\! +\!
\frac{\cosh\big(\sqrt{3\si_j^2-1}L\big)\th_j(x)}{\sinh\big(\sqrt{3\si_j^2\!-\!1}L\big)\!-\! e^{i\k_j^{(1)}\si_j L}\th_j(L)}\\
\ns&\ds \q +
\frac{\cosh\!\big(\!\sqrt{3\si_j^2\!-\!1}x\big)\sinh\!\big(\!\sqrt{3\si_j^2\!-\!1}L\big)\!-\sinh\!\big(\!\sqrt{3\si_j^2\!-\!1}x\big)\cosh\!\big(\!\sqrt{3\si_j^2\!-\!1}L\big)}{\sinh\big(\sqrt{3\si_j^2-1}L\big)-
e^{i\k_j^{(1)}\si_j L}\th_j(L)}
 e^{-i\k_j^{(1)}\si_j x}\\
 \ns&\ds \q -
\frac{\cosh\!\big(\!\sqrt{3\si_j^2\!-\!1}x\big)e^{i\k_j^{(1)}\si_j
(L-x)}\th_j(L) +
e^{i(\k_j^{(1)}+2\k_j^{(2)})\si_j
L}\th_j(x)}{\sinh\big(\sqrt{3\si_j^2-1}L\big)-e^{i\k_j^{(1)}\si_j
L}\th_j(L)}\\
\ns&\ds \q -
\frac{e^{ik_j^{(1)}\si_j(L-y)}\th_j(L)\big(\sinh\big(\sqrt{3\si_j^2-1}x\big)-e^{i\k_j^{(1)}\si_j
y}\th_j(x)
\big)}{\sinh\big(\sqrt{3\si_j^2-1}L\big)-e^{i\k_j^{(1)}\si_j
L}\th_j(L)} - e^{2i\k_j^{(2)}\si_j
x}\Bigg\}.
 \end{array}
\end{equation}
Then, using also \eqref{6.29-eq1}, we get
that
\begin{equation}\label{5.15-eq3.1}
\!\!\begin{array}{ll}\ds
\Bigg(\!\!\int_0^L\!\!
\Bigg|\hat\psi_j(x)\! -\!\a_{-j} \bigg[
\frac{2e^{i(\k_j^{(1)} +2\k_j^{(2)})\si_j
L} \sinh
\!\big(\!\sqrt{3\si_j^2\!-\!1}x\big)+
e^{\sqrt{3\si_j^2-1}(L-x)}}{2\sinh\!\big(\!\sqrt{3\si_j^2\!-\!1}L\big)-
2e^{i\k_j^{(1)}\si_j
L}\th_j(L)}e^{-i\k_j^{(1)}\si_j
x}\!\! - e^{2i\k_j^{(2)}\si_j x}\Bigg] \Bigg|^2\!dx\!\Bigg)^{\frac{1}{2}}\\
\ns \ds =  O\(\frac{1}{j^3}\) \qq \mbox{ as
} |j|\to +\infty.
\end{array}
\end{equation}

Now, we prove that
\begin{equation}\label{4.18-eq4}
\sum_{j\in\dbZ\setminus\{0\}}
\big|\hat\psi_{-j} - \f_j
\big|_{L^2(0,L)}^2 <\infty.
\end{equation}
Utilizing \eqref{6.20-eq4} and
straightforward computations, one can show
that
\begin{equation}\label{10.31-eq5}
\begin{array}{ll}\ds
\q\Bigg|\a_j \Bigg[ \frac{2e^{3i\tau_j
L}\sinh\big(\sqrt{3\tau_j^2-1}x\big) +
e^{\sqrt{3\tau_j^2-1}(L-x)}}{2\sinh\big(\sqrt{3\tau_j^2-1}L\big)}e^{-i\tau_j
x} - e^{2i\tau_j x}\Bigg]
\\
\ns\ds \q- \a_j \Bigg[
\frac{2e^{i(\k_j^{(1)}+2\k_j^{(2)})\si_{-j}
L}\sinh\big(\sqrt{3\si_{-j}^2-1}x\big)+ \!
e^{\sqrt{3\si_{-j}^2-1}(L-x)}}{2\sinh\big(\sqrt{3\si_{-j}^2-1}L\big)-
2e^{i\k_j^{(1)}\si_{-j}
L}\th_j(L)}e^{-i\k_j^{(1)}\si_{-j} x} -
e^{2i\k_j^{(2)}\si_{-j} x}\Bigg]\Bigg|
\\
\ns\ds \leq C\a_j \bigg(
\big|1-e^{(\sqrt{3\tau_j^2-1}-\sqrt{3\si_{-j}^2-1})x}\big|
+
\big|1-e^{(\k_j^{(1)}\si_{-j}-\tau_j)x}\big|
+
\big|1-e^{(\k_j^{(2)}\si_{-j}-\tau_j)x}\big|
+ e^{-\sqrt{3 \si_{-j}^2-1}L}  \bigg).
\end{array}
\end{equation}
Hence, using  \eqref{5.15-eq1},
\eqref{6.17-eq2} and \eqref{5.15-eq3.1}, we
see that there exists $C >0$ such that for
every $j \in\dbZ\setminus\{0\}$,
\begin{equation}\label{5.15-eq5}
\begin{array}{ll}\ds
\q \int_{0}^L |\hat \psi_{-j} -  \f_j|^2 dx \\
\ns\ds \leq C
\bigg(\Big|\sqrt{3\tau_j^2-1}-\sqrt{3\si_{-j}^2-1}\Big|^2
+ \big|\tau_j-\k_j^{(1)}\si_{-j}\big|^2  +
\big|\tau_j-\k_j^{(2)}\si_{-j}\big|^2 +
e^{-\sqrt{3 \si_{-j}^2-1}L} + \frac{1}{j^6}
\bigg)
\\
\ns \ds \leq C \frac{1}{j^4}.
\end{array}
\end{equation}

Since
$$
\sum_{j=1}^{+\infty} \frac{1}{j^4}
<+\infty,
$$
we have \eqref{4.18-eq4}.

Let $T>0$. Consider the following  control
system:
\begin{equation}\label{10.4-eq1}
\left\{
\begin{array}{ll}\ds
\vartheta_t + \vartheta_{xxx} + \vartheta_x = 0 &\mbox{ in } (0,T)\times(0,L),\\
\ns\ds \vartheta(t,0)=\vartheta(t,L)=0 &\mbox{ in } (0,T),\\
\ns\ds
\vartheta_x(t,L)-\vartheta_x(t,0)=\eta(t)
&\mbox{ in } (0,T),
\end{array}
\right.
\end{equation}
where $\eta(\cd)\in L^2(0,T)$ is the
control. Let $\tilde
\vartheta=A^{-1}\vartheta$ (recall
\eqref{DA} for the definition of $A$).
Then, we know $\tilde \vartheta$ solves
\begin{equation}\label{10.4-eq2}
\pa_t \tilde \vartheta = A\tilde \vartheta
- \eta(t)b,
\end{equation}
where $b(\cd)$ is the solution to
\begin{equation}\label{10.31-eq1}
\left\{
\begin{array}{ll}\ds
b_{xxx} + b_x =0 &\mbox{ in }(0,L),\\
\ns\ds b(0)=b(L)=0,\;b_x(L)-b_x(0)=1.
\end{array}
\right.
\end{equation}
Clearly, $b\in L^2(0,L)$. Let
$b=\sum_{j\in\dbZ\setminus\{0\}}b_j\f_j$.
Since the system \eqref{10.4-eq1} is
exactly controllable in $L^2(0,L)$ (see
\cite[Theorem 1.2]{Rosier1}), we get that
\eqref{10.4-eq1} is also exactly
controllable in $A^{-1}\big(L^2(0,L)\big)$.
In particular, $b_j\neq 0$ for every
$j\in\dbZ\setminus\{0\}$.

 For $j\in\dbZ\setminus\{0\}$, let $\tilde\psi_j\in D(A)$  be the solution of
\begin{equation}\label{10.4-eq3}
A\tilde\psi_j +
(-\lambda+i\mu_j)\tilde\psi_j =
-\frac{\a_{-j}}{\b_j}(-\lambda+i\mu_j) b.
\end{equation}
Then, for every $j\in\dbZ\setminus\{0\}$,
we have
\begin{equation}\label{10.4-eq3.1}
A^{-1}\tilde\psi_j =
-(-\lambda+i\mu_j)^{-1}\tilde\psi_j -
\frac{\a_{-j}}{\b_j}A^{-1}b
\end{equation}
and
\begin{equation}\label{10.4-eq10}
\hat\psi_j = \tilde \psi_j +
\frac{\a_{-j}}{\b_j}b.
\end{equation}

Assume that there exists
$\{a_j\}_{j\in\dbZ\setminus\{0\}}\in
\ell^2(\dbZ)$ such that
\begin{equation}\label{10.4-eq5}
\sum_{j\in\dbZ\setminus\{0\}}a_j\hat\psi_j=0.
\end{equation}
From \eqref{10.4-eq10} and
\eqref{10.4-eq5}, we obtain that
\begin{equation}\label{10.4-eq5.1}
\sum_{j\in\dbZ\setminus\{0\}}a_j\(\tilde
\psi_j + \frac{\a_{-j}}{\b_j} b\)=0.
\end{equation}

Applying $A^{-1}$ to \eqref{10.4-eq5.1},
and using \eqref{10.4-eq3.1}, one gets that
\begin{equation}\label{10.4-eq6}
\sum_{j\in\dbZ\setminus\{0\}}
a_j(-\lambda+i\mu_j)^{-1}\tilde\psi_j=0,
\end{equation}
which, together with \eqref{10.4-eq3.1},
implies that
\begin{equation}\label{10.4-eq6.1}
\bigg[\sum_{j\in\dbZ\setminus\{0\}}
a_j\frac{\a_{-j}}{\b_j}(-\lambda+i\mu_j)^{-1}\bigg]b-\sum_{j\in\dbZ\setminus\{0\}}
a_j(-\lambda+i\mu_j)^{-1}\hat\psi_j=0.
\end{equation}
Applying $A^{-1}$ to \eqref{10.4-eq6} and
using \eqref{10.4-eq3.1} again, we get
\begin{equation}\label{10.4-eq7}
\bigg[\sum_{j\in\dbZ\setminus\{0\}}
a_j\frac{\a_{-j}}{\b_j}(-\lambda+i\mu_j)^{-1}\bigg]A^{-1}b
+ \sum_{j\in\dbZ\setminus\{0\}}
a_j(-\lambda+i\mu_j)^{-2}\tilde\psi_j=0,
\end{equation}
which indicates that
\begin{equation}\label{10.4-eq7.1}
\bigg[\sum_{j\in\dbZ\setminus\{0\}}
a_j\frac{\a_{-j}}{\b_j}(-\lambda+i\mu_j)^{-1}\bigg]A^{-1}b
- \bigg[\sum_{j\in\dbZ\setminus\{0\}}
a_j\frac{\a_{-j}}{\b_j}(-\lambda+i\mu_j)^{-2}\bigg]
b + \sum_{j\in\dbZ\setminus\{0\}}
a_j(-\lambda+i\mu_j)^{-2}\hat\psi_j=0.
\end{equation}

By induction, one gets that, for every
positive integer $p$,
\begin{equation}\label{eqp}
\begin{array}{ll}\ds
\bigg(\sum_{j\in\dbZ\setminus\{0\}}
a_j\frac{\a_{-j}}{\b_j}(-\lambda+i\mu_j)^{-1}
\bigg)A^{-p}b +
\sum_{k=2}^{p}\bigg(\sum_{j\in\dbZ\setminus\{0\}}
a_j\frac{\a_{-j}}{\b_j}(-\lambda+i\mu_j)^{-k}(-1)^{p+1-k}
\bigg)A^{p+1-k}b \\
\ns\ds -  \sum_{j\in\dbZ\setminus\{0\}}
a_i\frac{\a_{-j}}{\b_j}(-\lambda+i\mu_j)^{-p-1}(-1)^{p+1}b
+
 \sum_{j\in\dbZ\setminus\{0\}}
a_i(-\lambda+i\mu_j)^{-p-1}(-1)^{p+1}\hat\psi_j=0.
\end{array}
\end{equation}
If
\begin{gather}
\label{premieresomme}
\sum_{j\in\dbZ\setminus\{0\}}
a_j\frac{\a_{-j}}{\b_j}(-\lambda+i\mu_j)^{-1}\not
=0,
\end{gather}
we get from \eqref{eqp} that
\begin{gather}
\label{inclusionpsi}
\{A^{-p}b\}_{p\in(\{0\}\cup\dbZ^+)}\subset\spa\{\hat\psi_j\}_{j\in\dbZ\setminus\{0\}}.
\end{gather}
If
$\overline{\spa\{\hat\psi_j\}_{j\in\dbZ\setminus\{0\}}}\neq
L^2(0,L)$, then we can find a nonzero
function
\begin{gather}
\label{dnot0}
d=\sum_{j\in\dbZ\setminus\{0\}}d_j\f_j\in
L^2(0,L)
\end{gather}
 such that
\begin{gather}
\label{scalarproduct0}
(h,d)_{L^2(0,L)}=0,\;\mbox{ for every }h\in
\spa\{\hat\psi_j\}_{j\in\dbZ\setminus\{0\}}.
\end{gather}
 From \eqref{inclusionpsi} and \eqref{scalarproduct0}, we obtain that $(A^{-p}b,d)_{L^2(0,L)}=0$
for every $p\in\{0\}\cup\dbZ^+$. Therefore,
we get
\begin{equation}\label{10.4-eq8}
\sum_{j\in\dbZ\setminus\{0\}}b_j(i\mu_j)^{-p}d_j=0
\mbox{ for all } p\in\{0\}\cup\dbZ^+.
\end{equation}
Let us define a complex variable function
$G(\cd):\mathbb{C}\to \mathbb{C}$ as
\begin{gather}
\label{defGcomplex}
G(z)=\sum_{j\in\dbZ\setminus\{0\}}d_jb_je^{(i\mu_j)^{-1}z},\q
z\in\dbC.
\end{gather}
Then, it is clear that $G(\cd)$ is a
holomorphic function. From
\eqref{10.4-eq8}, we see that
$$
G^{(p)}(0)=0 \;\mbox{ for every
}p\in\{0\}\cup\dbZ^+.
$$
Thus, we find that $G(\cd)= 0$, which
educes that $d_jb_j=0$ for all integer $j$.
Since $b_j\neq 0$ for every
$j\in\dbZ\setminus\{0\}$, we know that
$d_j=0$ for all integer $j$. Therefore, we
get $d=0$, which leads to a contradiction
with that $d\neq 0$. Hence,
\eqref{premieresomme} implies that
\begin{gather}
\label{density}
\overline{\spa\{\hat\psi_j\}_{j\in\dbZ\setminus\{0\}}}=L^2(0,L).
\end{gather}
If
$$\sum_{j\in\dbZ\setminus\{0\}}
a_j\frac{\a_{-j}}{\b_j}(-\lambda+i\mu_j)^{-1}=0
$$
and
$$
\sum_{j\in\dbZ\setminus\{0\}}
a_j\frac{\a_{-j}}{\b_j}(-\lambda+i\mu_j)^{-2}\neq
0,
$$
then by using \eqref{eqp} again, we obtain
that
$$
\{A^{-p}b\}_{p\in
\{0\}\cup\dbZ^+}\subset\spa\{\hat\psi_j\}_{j\in\dbZ\setminus\{0\}}.
$$
Then, by a similar argument, we find that
\eqref{density} again holds. Similarly, we
can get that, if there is a $p\in\dbZ^+$
such that
\begin{gather}\sum_{j\in\dbZ\setminus\{0\}} a_j
\frac{\a_{-j}}{\b_j}(-\lambda+i\mu_j)^{-p}
\neq 0,
\end{gather}
then \eqref{density} holds.

On the other hand, if
\begin{equation}\label{10.4-eq9}
\sum_{j\in\dbZ\setminus\{0\}} a_j
\frac{\a_{-j}}{\b_j}(-\lambda+i\mu_j)^{-p}
= 0 \mbox{ for every }p\in\dbZ^+,
\end{equation}
we  define a function
$$
\wt
G(z)\=\sum_{j\in\dbZ\setminus\{0\}}a_j\frac{\a_{-j}}{\b_j}(-\lambda+i\mu_j)^{-1}
e^{(-\lambda+i\mu_j)^{-1}z} \mbox{ for
every }z\in\dbC,
$$
and it is clear that $\wt G(\cd)$ is a
holomorphic function and
$$
\wt G^{(p)}(0)=0 \mbox{ for every }
p\in\{0\}\cup\dbZ^+,
$$
which implies that $\wt G(\cd)= 0$.
Therefore, we conclude that $a_j=0$ for
every $j\in\dbZ\setminus\{0\}$.

By the above argument, we know that either
$\{\hat\psi_j\}_{j\in\dbZ\setminus\{0\}}$
is $\o$-independent or is complete in
$L^2(0,L)$. Before proceeding our proof, we
recall the following two known results.

\begin{lemma}\cite[Page 45, Theorem 15]{1980-Young-book}\label{lm1}
Let $H$ be a separable Hilbert space and
let $\{e_j\}_{j\in\dbZ^+}$ be an
orthonormal basis for $H$. If
$\{f_j\}_{j\in\dbZ^+}$ is an
$\o$-independent sequence such that
$$
\sum_{j\in\dbZ^+}|f_j-e_j|_H^2<+\infty,
$$
then $\{f_j\}_{j\in\dbZ^+}$ is a Riesz
basis for $H$.
\end{lemma}
\begin{lemma}\cite[Page 40, Theorem 12]{1980-Young-book}\label{lm2}
Let $\{e_j\}_{j\in\dbZ^+}$ be a basis of a
Banach space $X$ and let
$\{f_j\}_{j\in\dbZ^+}$ be the associated
sequence of coefficient functionals. If
$\{b_j\}_{j\in\dbZ^+}$ is complete in $X$
and if
$$
\sum_{j\in\dbZ^+}|e_j-b_j|_{X}|f_j|_{X'}<
+\infty,
$$
then $\{b_j\}_{j\in\dbZ^+}$ is a basis for
$X$ which is equivalent to
$\{e_j\}_{j\in\dbZ^+}$.
\end{lemma}

Now we continue the proof. We first deal
with the case where
$\{\hat\psi_j\}_{j\in\dbZ\setminus\{0\}}$
is $\o$-independent. Let us take
$H=L^2(0,L)$ and put
$$
\left\{
\begin{array}{ll}
\ds e_{2j-1} = \f_j \mbox{ for } j\in\dbZ^+,\\
\ns\ds e_{2j} = \f_{-j} \mbox{ for }
j\in\dbZ^+,
\end{array}
\right. \qq \left\{
\begin{array}{ll}\ds
f_{2j-1} = \hat\psi_{-j} \mbox{ for } j\in\dbZ^+,\\
\ns\ds f_{2j} = \hat\psi_{j} \mbox{ for }
j\in\dbZ^+.
\end{array}
\right.
$$
Then, by \eqref{5.15-eq5}, the conditions
of Lemma \ref{lm1} are fulfilled. Thus, we
get that
$\{\hat\psi_j\}_{j\in\dbZ\setminus\{0\}}$
is a Riesz basis of $L^2(0,L)$.

Next, we consider the case that
$\{\hat\psi_j\}_{j\in\dbZ\setminus\{0\}}$
is complete in $L^2(0,L)$. Let us set
$X=L^2(0,L)$,
$$
\left\{
\begin{array}{ll}\ds e_{2j-1} = \f_j \mbox{ for } j\in\dbZ^+,\\
\ns\ds e_{2j} = \f_{-j} \mbox{ for }
j\in\dbZ^+,
\end{array}
\right. \qq \left\{
\begin{array}{ll}\ds f_{2j-1} = \f_j \mbox{ for } j\in\dbZ^+,\\
\ns\ds f_{2j} = \f_{-j} \mbox{ for }
j\in\dbZ^+,
\end{array}
\right. \qq \left\{
\begin{array}{ll}\ds b_{2j-1} = \hat\psi_{-j} \mbox{ for } j\in\dbZ^+,\\
\ns\ds b_{2j} = \hat\psi_{j} \mbox{ for }
j\in\dbZ^+.
\end{array}
\right.
$$
Then, using \eqref{5.15-eq5} once more, it
is easy to see that the conditions of Lemma
\ref{lm2} are fulfilled. Therefore, we get
that
$\{\hat\psi_j\}_{j\in\dbZ\setminus\{0\}}$
is a Riesz basis of $L^2(0,L)$.

Now we give an estimate of
$\{c_j\}_{j\in\dbZ\setminus\{0\}}$. From
\eqref{estimatetauj} and \eqref{8.30-eq15},
we get that
$$
\frac{1}{C} j^2\leq \int_0^L
|\f_{j,y}(y)|^2 dy \leq C j^2.
$$
This, together with the fact that
$\{\hat\psi_j\}_{j\in\dbZ\setminus\{0\}}$
is a Riesz basis of $L^2(0,L)$ and
$k_3(\cd,\cd)\in H_0^1((0,L)\times (0,L))$,
implies that
$$
\begin{array}{ll}\ds
+\infty\3n&\ds> \int_0^L\int_0^L
|k_{3,y}(x,y)|^2dxdy \geq \int_0^L\int_0^L
\|\sum_{j\in\dbZ\setminus\{0\}}
c_j\frac{\b_j}{\a_{-j}}\hat\psi_j(x)\f_{j,y}(y)\|^2dxdy\\
\ns&\ds \geq C\int_0^L
\sum_{j\in\dbZ\setminus\{0\}}\|
c_j\frac{\b_j}{\a_{-j}} \f_{j,y}(y)\|^2
dy\geq C\sum_{j\in\dbZ\setminus\{0\}}
\Big|c_j\frac{\b_j}{\a_{-j}}j\Big|^2.
\end{array}
$$
Hence, we find that
\begin{equation}\label{10.5-eq7}
\Bigg\{c_j\frac{\b_j}{\a_{-j}}j\Bigg\}_{j\in\dbZ\setminus\{0\}}\in
\ell^2(\dbZ\setminus\{0\}).
\end{equation}
From \eqref{fjx(0)not0},
\eqref{estimatetauj}, \eqref{6.29-eq1} and
\eqref{8.30-eq5}, we get the existence of
$C>0$ such that
\begin{equation}
\label{asympfjy-better} C^{-1}|j|\leq
|\f_{j,y}(0)| \leq C |j|, \q \forall j\in
\mathbb{Z}\setminus\{0\}.
\end{equation}
From \eqref{10.5-eq7} and
\eqref{asympfjy-better}, we find that
\begin{equation}\label{10.5-eq8}
\Bigg\{c_j\frac{\b_j}{\a_{-j}}\f_{j,y}(0)\Bigg\}_{j\in\dbZ\setminus\{0\}}\in
\ell^2(\dbZ\setminus\{0\}).
\end{equation}
Using \eqref{k3=}, $k_{3,y}(x,0)=0$ and the
second inequality of
\eqref{asympfjy-better}, we find that
$$
\sum_{j\in\dbZ\setminus\{0\} }
\psi_j(\cd)\f_{j,y}(0)=\sum_{j\in\dbZ\setminus\{0\}
}
c_j\frac{\b_j}{\a_{-j}}\f_{j,y}(0)\hat\psi_j(\cd)=0
\q\mbox{ in } L^2(0,L).
$$
Using that $\{\hat\psi_j\}_{j\in \dbZ }$ is
a Riesz basis of $L^2(0,L)$, we get
$$
c_j\frac{\b_j}{\a_{-j}}\f_{j,y}(0)=0 \mbox{
for every } j\in\dbZ\setminus\{0\}.
$$
 Then, we see $\f_{j,y}(0)=0$ if $c_j\neq 0$.
However, by \cite[Lemma 3.5]{Rosier1}, this
is impossible since $L\notin\cN$. Thus, we
get that $c_j=0$ for every
$j\in\dbZ\setminus\{0\}$, which implies
that the equation \eqref{system4} admits a
unique solution $k_3=0$. Therefore, we
obtain that the equation \eqref{system3}
admits at most one solution. This concludes
Step 1.

\vspace{0.3cm}

\begin{center}
\textbf{Step 2: proof of the existence of a
solution to \eqref{system3} with the
required regularity}
\end{center}

Denote by $D(A)'$ the dual space of $D(A)$
with respect to the pivot space $L^2(0,L)$.
Let
$h(\cd)=\sum_{j\in\dbZ\setminus\{0\}}h_j\f_j(\cd)\in
D(A)$, i.e., $|h(\cd)|_{D(A)}^2 =
\sum_{j\in\dbZ\setminus\{0\}}|h_j\mu_j|^2<+\infty$,
we have that
$$
\sum_{j\in\dbZ\setminus\{0\}}|h_j\f_{j}'(0)|
\leq
\(\sum_{j\in\dbZ\setminus\{0\}}|h_j\mu_j|^2\)^\frac{1}{2}\(\sum_{j\in\dbZ\setminus\{0\}}|\f_{j}'(0)\mu_j^{-1}|^2\)^\frac{1}{2}<+\infty.
$$
Hence
\begin{equation}
\sum_{j=-l_1}^{l_2}\f_{j}'(0)\ov{\f_j(\cd)}
\text{ is converging in $D(A)'$ as $l_1$
and $l_2$ tend to $+\infty$},
\end{equation}
which  allows us to define
\begin{equation}
\label{defhatb} \hat b(\cd)
\=\sum_{j\in\dbZ\setminus\{0\}}\f_{j}'(0)\ov{\f_j(\cd)}
\in D(A)'.
\end{equation}
Further, it is clear that
$$
(h,\hat b)_{D(A),D(A)'} =
\sum_{j\in\dbZ\setminus\{0\}}h_j\f_{j}'(0)=h'(0)=-\d'_0(h),
$$
which implies that
\begin{equation}\label{11.8-eq5}
\hat b = -\d'_0 \q\mbox{ in } D(A)'.
\end{equation}

Let
\begin{equation}\label{11.7-eq1}
a_{j} \= \frac{\l}{\f'_{j}(0)},\;\;\phi_j
\= \bar\f_{j} +
(-A-i\mu_{j}+\l)^{-1}(a_{j}\hat b - \l
\bar\f_{j}) \;\mbox{ for }
j\in\dbZ\setminus\{0\}.
\end{equation}
From \eqref{5.22-eq1},
\eqref{asympfjy-better} and
\eqref{11.7-eq1}, we have that
\begin{equation}\label{10.30-eq1}
\begin{array}{ll}\ds
\q\int_{0}^L \Big|\ov{\f_{j}(x)} -
\phi_{j}(x)\Big|^2 dx\\
\ns\ds = \int_{0}^L
\big|(-A-i\mu_{j}+\l)^{-1}(a_{j}\hat b -
\l \bar\f_{j})\big|^2 dx\\
\ns\ds =  \int_{0}^L
\Big|\sum_{k\in\dbZ\setminus\{0,j\}}(i\mu_k
-
i\mu_{j}+\l)^{-1}\frac{\l\f_k'(0)}{\f_{j}'(0)}\bar
\f_k(x)\Big|^2dx
\\
\ns\ds = \sum_{k\in\dbZ\setminus\{0,j\}}
\Big|(i\mu_k-i\mu_{j}+\l)^{-1}\frac{\l\f_k'(0)}{\f_{j}'(0)} \Big|^2\\
\ns\ds\leq
C\l^2\sum_{k\in\dbZ\setminus\{0,j\}}
\left|k^3-j^3 \right|^{-2}j^{-2}k^2
\q\mbox{ for all }j\in\dbZ\setminus\{0\}.
\end{array}
\end{equation}
From \eqref{10.30-eq1}, we obtain, for
$j>0$,
\begin{equation}\label{11.7-eq2}
\begin{array}{ll}\ds
\q\int_{0}^L \Big|\ov{\f_{j}(x)} -
\phi_{j}(x)\Big|^2 dx
\\
\ns\ds \leq
C\bigg(\sum_{k>2j}\frac{k^2}{\big|k^3-j^3\big|^{2}j^{2}}
+ \sum_{j< k\leq
2j}\frac{k^2}{\big|k^3-j^3\big|^{2}j^{2}} +
\sum_{-j< k <j,\,k\neq 0}\frac{k^2}{\big|k^3-j^3\big|^{2}j^{2}}\\
\ns\ds \q + \sum_{k\leq-
j}\frac{k^2}{\big|k^3-j^3\big|^{2}j^{2}}\bigg).
\end{array}
\end{equation}
Now we estimate the terms in the right hand
side of \eqref{11.7-eq2}. First, we have
that
\begin{equation}\label{11.7-eq3}
\sum_{k>2j}\frac{k^2}{\big|k^3-j^3\big|^{2}j^{2}}\leq
\sum_{k>2j}\frac{8}{7k^{4}j^{2}} \leq
\frac{8}{7j^4}\sum_{k>2j}\frac{1}{k^{2}}\leq
\frac{1}{j^5}.
\end{equation}
Next,
\begin{equation}\label{11.7-eq4}
\begin{array}{ll}\ds
\sum_{j< k\leq
2j}\frac{k^2}{\big|k^3-j^3\big|^{2}j^{2}}\leq
4\sum_{j< k\leq
2j}\frac{1}{\big|k^3-j^3\big|^{2}} =
4\sum_{1\leq l\leq
j}\frac{1}{\big|(j+l)^3-j^3\big|^{2}}\\
\ns\ds \leq 4\sum_{1\leq l\leq
j}\frac{1}{\big|3j^2l+3jl^2+l^3\big|^{2}}\leq
\frac{4}{9}\sum_{-j\leq l\leq
-1}\frac{1}{j^2}\frac{1}{(j+l)^2l^2}=
\frac{4}{9}\sum_{1\leq l\leq
j}\frac{1}{j^2}\[\frac{1}{j}\(\frac{1}{l}-\frac{1}{j+l}\)\]^2\\
\ns\ds \leq \frac{8}{9}\sum_{1\leq l\leq
j}\frac{1}{j^4}\(\frac{1}{l^2}+\frac{1}{(j-l)^2}\)\leq
\frac{8\pi^2}{27}\frac{1}{j^4}.
\end{array}
\end{equation}
Similarly, we can get that
\begin{equation}\label{11.7-eq5}
\sum_{-j<k<j,\, k\neq
0}\frac{k^2}{\big|k^3-j^3\big|^{2}j^{2}}
\leq \frac{8\pi^2}{27}\frac{1}{j^4}
\end{equation}
and
\begin{equation}\label{11.7-eq6}
\sum_{k\leq
-j}\frac{k^2}{\big|k^3+j^3\big|^{2}j^{2}}\leq
\frac{8}{j^5}.
\end{equation}
From \eqref{11.7-eq2} to \eqref{11.7-eq6},
we know that there is a constant $C>0$ such
that for all positive integer $j$, it holds
that
\begin{equation}\label{10.30-eq7}
\int_{0}^L \Big|\ov{\f_{j}(x)} -
\phi_{j}(x)\Big|^2 dx \leq \frac{C}{j^4}.
\end{equation}
Similarly, we can prove that there is a
constant $C>0$ such that for all negative
integer $j$, it holds that
\begin{equation}\label{10.30-eq8}
\int_{0}^L \Big|\ov{\f_{j}(x)} -
\phi_{j}(x)\Big|^2 dx \leq \frac{C}{j^4}.
\end{equation}
Further, by similar arguments, we can
obtain that
\begin{equation}\label{10.30-eq2}
\int_{0}^L \Big|\ov{\f_{j,x}(x)} -
 \phi_{j,x}(x)\Big|^2 dx \leq
\frac{C}{j^2}\q\mbox{ for all
}j\in\dbZ\setminus\{0\}.
\end{equation}
With the same strategy to prove that
$\{\hat\psi_j\}_{j\in\dbZ\setminus\{0\}}$
is a Riesz basis of $L^2(0,L)$, we also can
show that
\begin{equation}\label{10.31-eq3}
\Big\{-(\l-i\mu_{j})A^{-1}\phi_{j}\Big\}_{j\in\dbZ\setminus\{0\}}
\mbox{ is a Riesz basis   of }L^2(0,L).
\end{equation}
 From \eqref{5.22-eq1} and \eqref{asympfjy-better}, we  get that

$$\sum_{j\in\dbZ\setminus\{0\}}\f_{j}'(0)(i\mu_j)^{-1}\bar\f_j \text{ is converging in }
L^2(0,L).$$  From \eqref{10.31-eq3}, there
is $\{\hat
c_j\}_{j\in\dbZ\setminus\{0\}}\in
l^2(\dbZ\setminus\{0\})$ such that
$$
-\sum_{j\in\dbZ\setminus\{0\}} \hat c_j
(\l-i\mu_{j})A^{-1}\phi_{j}=\sum_{j\in\dbZ\setminus\{0\}}\f_{j,y}(0)(i\mu_j)^{-1}\bar\f_j
\q\mbox{ in }\; L^2(0,L),
$$
which, together with \eqref{defhatb} and
\eqref{11.8-eq5}, implies that
\begin{equation}\label{11.2-eq2}
-\sum_{j\in\dbZ\setminus\{0\}} \hat c_j
(\l-i\mu_{j})
\phi_{j}=\sum_{j\in\dbZ\setminus\{0\}}\f_{j,y}(0)\bar\f_j=-\d'_0
\q\mbox{ in }\; D(A)'.
\end{equation}
Since $\f_j\in D(A)$, from
\eqref{11.2-eq2}, we find that
\begin{equation}\label{11.2-eq1}
-\hat c_{j} (\l-i\mu_{j}) \int_0^L \f_j(x)
\phi_{j}(x)dx-\sum_{k\in\dbZ\setminus\{0,j\}}
\hat c_k (\l-i\mu_{k}) \int_0^L \f_j(x)
\phi_k(x)dx=\f_{j,y}(0).
\end{equation}
For $j\neq k$, we get, using
\eqref{defhatb} and \eqref{11.7-eq1},
$$
\begin{array}{ll}\ds
\q\hat c_k (\l-i\mu_{k})
\int_0^L \f_j(x) \phi_k(x)dx\\
\ns\ds =\hat c_k (\l-i\mu_{k})
\int_0^L \f_j(x)\big\{\ov{\f_{k}(x)}+(-A-i\mu_{k}+\l)^{-1}\big[a_{k}\hat b(x)-\l\ov{\f_{k}(x)}\,\big]\big\} dx\\
\ns\ds = \hat c_k (\l-i\mu_{k}) \int_0^L
\f_j(x)(-A-i\mu_{k}+\l)^{-1}a_{k} \f_j'(0)
\ov{\f_j(x)}dx\\
\ns\ds =  \hat c_k
(\l-i\mu_{k})(i\mu_j-i\mu_k+\l)^{-1}a_{k} \f_j'(0) \\
\ns\ds = \hat c_k
\l\frac{\l-i\mu_{k}}{i\mu_j-i\mu_k+\l}\frac{\f'_{j}(0)}{\f'_{k}(0)}.
\end{array}
$$
This, together with \eqref{11.2-eq1},
implies that
\begin{equation}\label{11.3-eq3}
-\hat c_{j} (\l-i\mu_{j})  -
\sum_{k\in\dbZ\setminus\{0,j\}}\hat c_k
\l\frac{\l-i\mu_{k}}{i\mu_j-i\mu_k+\l}\frac{\f'_{j}(0)}{\f'_{k}(0)}=\f'_{j}(0)\;
\mbox{ for all } j\in\dbZ\setminus\{0\}.
\end{equation}
For $j\in\dbZ\setminus\{0\}$ and
$k\in\dbZ\setminus\{0\}$, let
\begin{equation}\label{11.8-eq4}
c_j=-\frac{\hat
c_{j}(\l-i\mu_{j})}{\f'_{j}(0)}, \, a_{jk}=
\frac{1}{i\mu_j-i\mu_{k}+\l}.
\end{equation}
 From
\eqref{11.3-eq3} and \eqref{11.8-eq4}, we
get that
\begin{equation}\label{11.4-eq1}
c_{j}+\l\sum_{k\in\dbZ\setminus\{0,j\}}
a_{jk}c_k =1.
\end{equation}
Let us now estimate $c_j$.   From
\eqref{5.22-eq1}, \eqref{asympfjy-better}
and \eqref{11.4-eq1}, we have, for every
$j\in \dbZ\setminus\{0\}$,
\begin{equation}\label{11.7-eq15}
\begin{array}{ll}\ds
\Big|\sum_{k\in\dbZ\setminus\{0,j\}}
a_{jk}c_k\Big| \3n&\ds=
\Big|\sum_{k\in\dbZ\setminus\{0,j\}}\hat
c_k
\frac{\l-i\mu_{k}}{i\mu_j-i\mu_k+\l}\frac{1}{\f'_{k}(0)}\Big|\\
\ns&\ds \leq
C\sum_{k\in\dbZ\setminus\{0,j\}}\Big|\hat
c_k \frac{\l-ik^3}{j^3 -
k^3}\frac{1}{k}\Big|\\
\ns&\ds \leq
C\(\sum_{k\in\dbZ\setminus\{0,j\}}|\hat
c_k|^2\)^\frac{1}{2}\(\sum_{k\in\dbZ\setminus\{0,j\}}\Big|
\frac{\l-ik^3}{j^3 -
k^3}\frac{1}{k}\Big|^2\)^\frac{1}{2}.
\end{array}
\end{equation}
For positive integer $j$,
\begin{equation}\label{11.7-eq16}
\begin{array}{ll}\ds
\sum_{k\in\dbZ\setminus\{0,j\}}\Big|
\frac{\l-ik^3}{j^3 -
k^3}\frac{1}{k}\Big|^2\3n&\ds\leq \sum_{k>
2j}\Big| \frac{\l-ik^3}{j^3 -
k^3}\frac{1}{k}\Big|^2 + \sum_{j< k\leq
2j}\Big|
\frac{\l-ik^3}{j^3 - k^3}\frac{1}{k}\Big|^2\\
\ns&\ds\q + \sum_{-j< k < j,\,k\neq 0}\Big|
\frac{\l-ik^3}{j^3 - k^3}\frac{1}{k}\Big|^2
+ \sum_{k\leq -j}\Big| \frac{\l-ik^3}{j^3 -
k^3}\frac{1}{k}\Big|^2.
\end{array}
\end{equation}
We now estimate the four terms in the right
hand side of \eqref{11.7-eq16}. First, we
have that
\begin{equation}\label{11.7-eq17}
\sum_{k>2j}\Big| \frac{\l-ik^3}{j^3 -
k^3}\frac{1}{k}\Big|^2 \leq
\frac{8}{7}\sum_{k>2j}
\frac{1+\l^2}{k^2}\leq \frac{C}{j}.
\end{equation}
Second,
\begin{equation}\label{11.7-eq18}
\begin{array}{ll}\ds
\sum_{j< k\leq 2j}\Big| \frac{\l-ik^3}{j^3
- k^3}\frac{1}{k}\Big|^2 \leq  C+ \sum_{j<
k\leq 2j}
\Big|\frac{4j^2}{j^3-k^3}\Big|^2\leq C+ \sum_{1\leq l\leq j} \Big|\frac{4j^2}{j^3-(j+l)^3}\Big|^2\\
\ns\ds \leq  C + \sum_{1\leq l\leq j}
\Big|\frac{4j^2}{3j^2l + 3jl^2+l^3}\Big|^2
\leq C+ \frac{16}{9}\sum_{1\leq l\leq j}
\frac{1}{l^2} \leq C.
\end{array}
\end{equation}
With similar arguments, we can obtain that
\begin{equation}\label{11.7-eq19}
\sum_{-j< k< j,\,k\neq 0}\Big|
\frac{\l-ik^3}{j^3 - k^3}\frac{1}{k}\Big|^2
\leq C
\end{equation}
and
\begin{equation}\label{11.7-eq20}
\sum_{k \leq -j}\Big| \frac{\l-ik^3}{j^3 -
k^3}\frac{1}{k}\Big|^2 \leq C.
\end{equation}
From \eqref{11.7-eq15} to
\eqref{11.7-eq20}, we find that there is a
constant $C>0$ such that for all positive
integer $j$,
\begin{equation}\label{11.7-eq20.1}
\Big|\sum_{k\in\dbZ\setminus\{0,j\}}
a_{jk}c_k\Big|\leq C.
\end{equation}
Similarly, we can prove that for all
negative integer $j$,
\begin{equation}\label{11.7-eq20.2}
\Big|\sum_{k\in\dbZ\setminus\{0,j\}}
a_{jk}c_k\Big|\leq C.
\end{equation}
Combining \eqref{11.4-eq1},
\eqref{11.7-eq20.1} and
\eqref{11.7-eq20.2}, we get that there is a
constant $C>0$ such that for all
$j\in\dbZ\setminus\{0\}$,
\begin{equation}\label{11.7-eq21}
|c_j|\leq C.
\end{equation}

We now estimate $|c_j|$ for $|j|$ large.
From \eqref{5.22-eq1} and \eqref{11.8-eq4},
 we get  that for $j>0$,
\begin{equation}\label{11.7-eq7}
\begin{array}{ll}\ds
\sum_{k\in\dbZ\setminus\{0,j\}}\!|a_{jk}|\!
 =\!\!\sum_{k\in\dbZ\setminus\{0,j\}}\!\frac{1}{|i\mu_j\!-\!i\mu_k\!+\!\l|}\leq\!\!
\sum_{k\in\dbZ\setminus\{0,j\}}\frac{1}{|\mu_j-\mu_k
|}  \leq
\sum_{k\in\dbZ\setminus\{0,j\}}\frac{C}{|j^3-k^3|}\\
\ns\ds \leq
\sum_{k>2j}\frac{C}{\big|j^3-k^3\big|} +
\sum_{j< k\leq
2j}\frac{C}{\big|j^3-k^3\big|} +
\sum_{-j<k<j,\,k\neq
0}\frac{C}{\big|j^3-k^3\big|} +
\sum_{k\leq- j}\frac{C}{\big|j^3-k^3\big|}.
\end{array}
\end{equation}
We estimate the terms in the last line of
\eqref{11.7-eq7} one by one. First,
\begin{equation}\label{11.7-eq8}
\sum_{k>2j}\frac{1}{|j^3-k^3|}\leq
\frac{8}{7} \sum_{k>2j}\frac{1}{k^3}\leq
\frac{8}{7}
\sum_{k>2j}\frac{1}{k}\frac{1}{k(k-1)}\leq
\frac{8}{7}\frac{1}{(2j+1)2j} \leq
\frac{2}{7j^2}.
\end{equation}
Second,
\begin{equation}\label{11.7-eq9}
\begin{array}{ll}\ds
\sum_{j< k\leq
2j}\frac{1}{|j^3-k^3|}\3n&\ds\leq \sum_{1<
l\leq j}\frac{1}{|j^3-(j+l)^3|}\leq
\sum_{1< l\leq
j}\frac{1}{3j^2l+3jl^2+l^3}\\
\ns&\ds \leq \frac{1}{3}\sum_{1< l\leq
j}\frac{1}{j^2l+jl^2}\leq
\frac{1}{3j^2}\sum_{1< l\leq
j}\(\frac{1}{j+l} + \frac{1}{l}\)\\
\ns&\ds \leq \frac{2\ln|j|}{3j^2}.
\end{array}
\end{equation}
Similarly, we can obtain that
\begin{equation}\label{11.7-eq10}
\sum_{-j<k<j,\,k\neq
0}\frac{1}{\big|j^3-k^3\big|}\leq
\frac{4\ln|j|}{3j^2}
\end{equation}
and
\begin{equation}\label{11.7-eq11}
\sum_{k\leq-j}\frac{1}{\big|j^3-k^3\big|}\leq
\frac{1}{j^2}.
\end{equation}
From \eqref{11.7-eq7} to \eqref{11.7-eq11},
we know there is a constant $C>0$ such that
for all positive integer $j$,
\begin{equation}\label{11.7-eq12}
\sum_{k\in\dbZ\setminus\{0,j\}}|a_{jk}|
\leq \frac{C\ln|j|}{j^2}.
\end{equation}
By similar arguments, we also can show that
for all negative integer $j$,
\begin{equation}\label{11.7-eq13}
\sum_{k\in\dbZ\setminus\{0,j\}}|a_{jk}|
\leq \frac{C\ln|j|}{j^2}.
\end{equation}
Combining \eqref{11.4-eq1},
\eqref{11.7-eq21}, \eqref{11.7-eq12} and
\eqref{11.7-eq13}, we obtain that
\begin{equation}\label{11.7-eq22}
1-\frac{C\ln|j|}{j^2} \leq c_j \leq
1+\frac{C\ln|j|}{j^2}.
\end{equation}

We now turn to the definition of
$k(\cd,\cd)$. From \eqref{10.30-eq7} and
\eqref{11.7-eq22}, one has
\begin{equation}\label{11.7-eq23}
\begin{array}{lcl}\ds
\q\sum_{j\in\dbZ\setminus\{0\}}
\int_0^L\big|\ov{\f_{j}(x)}-c_j\phi_{j}(x)\big|^2
dx\\
\ns\ds \leq \ds
2\sum_{j\in\dbZ\setminus\{0\}}
\int_0^L\big|(1-c_j)\ov{\f_{j}(x)}\big|^2
dx  + 2\sum_{j\in\dbZ\setminus\{0\}}
\int_0^L|c_j|^2\big|\ov{\f_{j}(x)}-\phi_{j}(x)\big|^2
dx
\\
\ns\ds  \leq
C\sum_{j\in\dbZ\setminus\{0\}}\(\frac{\ln^2|j|}{j^4}
+\frac{1}{j^4}\) <+\infty.
\end{array}
\end{equation}
Inequality \eqref {11.7-eq23} allows us to
define $k(\cd,\cd)\in
L^2((0,L)\times(0,L))$ by
\begin{equation}
\label{defk}
k(x,y)\=\sum_{j\in\dbZ\setminus\{0\}}[\ov{\f_{j}(x)}-c_j\phi_{j}(x)]\f_j(y),\q
0\leq x,y\leq L.
\end{equation}
Let us prove that $k(\cd,\cd)\in
H^1((0,L)\times (0,L))$.
Thanks to \eqref{estimateLinftyfjx(x)},
\eqref{10.30-eq2}  and \eqref{11.7-eq22},
we find that
\begin{equation}\label{11.7-eq24}
\begin{array}{ll}\ds
\q\int_0^L\int_0^L\Big|\sum_{j\in\dbZ\setminus\{0\}}
[\ov{\f_{j}'(x)}-c_j\phi'_{j}(x)]\f_j(y)\Big|^2
dxdy\\
\ns\ds  =  \sum_{j\in\dbZ\setminus\{0\}}
\int_0^L\big|\ov{\f_{j}'(x)}-c_j\phi_{j}'(x)\big|^2
dx \\
\ns\ds \leq\sum_{j\in\dbZ\setminus\{0\}}
\int_0^L\big|(1-c_j)\ov{\f_{j}'(x)}\big|^2
dx + \sum_{j\in\dbZ\setminus\{0\}}
\int_0^L|c_j|^2\big|\ov{\f_{j}'(x)}-\phi_{j}'(x)\big|^2
dx
\\
\ns\ds \leq
C\sum_{j\in\dbZ\setminus\{0\}}\(\frac{\ln^2|j|}{j^2}
+\frac{1}{j^2}\) <+\infty,
\end{array}
\end{equation}
which shows that
\begin{equation}\label{kxL2}
  k_x(\cdot,\cdot)\in L^2((0,L)\times(0,L)).
\end{equation}

Utilizing \eqref{5.22-eq1},
\eqref{11.7-eq1} and  \eqref{11.7-eq21}, we
get
\begin{equation}\label{11.7-eq25}
\begin{array}{ll}\ds
\q\int_0^L\int_0^L\Big|\sum_{j\in\dbZ\setminus\{0\}}
c_j[\ov{\f_{j}(x)}-\phi_{j}(x)]\f_j'(y)\Big|^2dxdy
\\
\ns\ds  =\int_0^L
\int_0^L\sum_{j,k\in\dbZ\setminus\{0\}}
c_jc_k\big[\ov{\f_{j}(x)}-\phi_j(x)\big]\big[
\f_{k}(x)-\ov{\phi_k(x)}\,\big]\f_j'(y)\ov{\f_k'(y)}dxdy \\
\ns\ds =
\int_0^L\int_0^L\sum_{j,k\in\dbZ\setminus\{0\}}\!\!c_jc_k
\!(-A\!-\!i\mu_{j}\!+\!\l)^{-1}(a_{j}\hat b
\!-\! \l
\ov{\f_{j}})(x)\overline{(-A\!-\!i\mu_{k}\!+\!\l)^{-1}(a_{k}\hat
b \!-\!
\l \ov{\f_{k}})}(x)\\
\ns\ds \qq\qq\qq\q \times \f_j'(y)\ov{\f_k'(y)}dxdy  \\
\ns\ds = \int_0^L\int_0^L
\sum_{j,k\in\dbZ\setminus\{0\}}c_jc_k
\[\sum_{l\in\dbZ\setminus\{0,j\}}(i\mu_l-i\mu_{j}+\l)^{-1}
\frac{\l\f'_l(0)}{\f'_{j}(0)}\ov{\f_l(x)}\]\\
\ns\ds\q
\times\overline{\[\sum_{l\in\dbZ\setminus\{0,k\}}(i\mu_l-i\mu_{k}+\l)^{-1}
\frac{\l\f'_l(0)}{\f'_{k}(0)}\ov{\f_l(x)}\]}\f_j'(y)\ov{\f_k'(y)}dxdy \\
\ns\ds \leq
C\l^2\!\!\sum_{j,k\in\dbZ\setminus\{0\}}
\sum_{l\in\dbZ\setminus\{0,j,k\}}\big|(i\mu_l-i\mu_{j}+\l)^{-1}
\f'_l(0) (i\mu_l-i\mu_{k}+\l)^{-1}
\f'_l(0)\big|\\
\ns\ds \q \times\int_0^L\big|\f_j'(y)\f'_j(0)^{-1}\f_k'(y)\f'_k(0)^{-1}\Big|dy\\
\ns\ds \leq
C\l^2\sum_{l\in\dbZ\setminus\{0\}}
\sum_{j\in\dbZ\setminus\{0,l\}}\big|(i\mu_l-i\mu_{j}+\l)^{-1}
\f'_l(0)\big|
\sum_{k\in\dbZ\setminus\{0,l\}}\big|(i\mu_l-i\mu_{k}+\l)^{-1}
\f'_l(0)\big|.
\end{array}
\end{equation}
Similar to the proof of \eqref{11.7-eq12},
we can obtain that
\begin{equation}\label{11.7-eq26}
\left\{
\begin{array}{ll}\ds
\sum_{j\in\dbZ\setminus\{0,l\}}\big|(i\mu_l-i\mu_{j}+\l)^{-1}
\f'_l(0)\big| \leq \frac{C\ln|l||\f'_l(0)|}{l^2}\leq \frac{C\ln|l|}{|l|},\\
\ns\ds
\sum_{k\in\dbZ\setminus\{0,l\}}\big|(i\mu_l-i\mu_{k}+\l)^{-1}
\f'_l(0)\big|  \leq
\frac{C\ln|l||\f'_l(0)|}{l^2}\leq
\frac{C\ln|l|}{|l|}.
\end{array}
\right.
\end{equation}
Combining \eqref{11.7-eq25} and
\eqref{11.7-eq26}, one has that
\begin{equation}\label{11.7-eq27}
\begin{array}{ll}\ds
\int_0^L\int_0^L\Big|\sum_{j\in\dbZ\setminus\{0\}}
c_j\big[\,\ov{\f_{j}(x)}-\phi_{j}(x)\big]\f_j'(y)\Big|^2dxdy
\leq
C\l^2\!\!\sum_{l\in\dbZ\setminus\{0\}}\!\frac{\ln^2|l|}{l^2}<+\infty.
\end{array}
\end{equation}
 From \eqref{estimateLinftyfjx(x)}, \eqref{11.7-eq22} and  \eqref{11.7-eq27}
\begin{equation}\label{11.8-eq1}
\begin{array}{ll}\ds
\q\int_0^L\int_0^L\Big|\sum_{j\in\dbZ\setminus\{0\}}
\big[\,\ov{\f_{j}(x)}-c_j\phi_{j}(x)\big]\f_j'(y)\Big|^2
dxdy\\
\ns\ds  \leq
2\int_0^L\int_0^L\Big|\sum_{j\in\dbZ\setminus\{0\}}(1-c_j)\ov{\f_{j}(x)}\f_j'(y)\Big|^2
dxdy\\
\ns\ds \q +
2\int_0^L\int_0^L\Big|\sum_{j\in\dbZ\setminus\{0\}}
c_j\big[\,\ov{\f_{j}(x)}-\phi_{j}(x)\big]\f_j'(y)\Big|^2dxdy
\\
\ns\ds \leq 2\sum_{j\in\dbZ\setminus\{0\}}
\int_0^L\big|(1-c_j)\f_{j}'(y)\big|^2 dy +
2\int_0^L\int_0^L\Big|\sum_{j\in\dbZ\setminus\{0\}}
c_j\big[\,\ov{\f_{j}(x)}-\phi_{j}(x)\big]\f_j'(y)\Big|^2dxdy
\\
\ns\ds \leq
C\sum_{j\in\dbZ\setminus\{0\}}\frac{\ln^2|j|}{j^2}
 <+\infty.
\end{array}
\end{equation}
From \eqref{11.7-eq23}, \eqref{11.7-eq24}
and \eqref{11.8-eq1}, we see that
\begin{equation}
\label{kyL^2} k_y(\cd,\cd)\in
L^2((0,L)\times (0,L)),
\end{equation}
which, together with \eqref{kxL2}, shows
that $k(\cd,\cd)\in H^1((0,L)\times
(0,L))$.  Clearly, $k(\cd,\cd)=0$ on the
boundary of $(0,L)\times (0,L)$. Thus, we
conclude that $k(\cd,\cd)\in
H^1_0((0,L)\times (0,L))$.

Let us define
\begin{equation}
\label{defkn} k^{(n)}(x,y) \=
\sum_{0<|j|\leq
n}\big[\,\ov{\f_{j}(x)}-c_j\phi_j(x)\big]\f_j(y)
\;\mbox{ in }(0,L)\times (0,L).
\end{equation}
Simple estimates show that
\begin{equation}\label{10.13-eq4}
\left(x\in (0,L)\mapsto
k^{(n)}_{x}(x,\cd)\in L^2(0,L)\right)\text{
is in } C^0([0,L]; L^2(0,L)).
\end{equation}

For any $m,n\in\dbZ^+$, $m<n$, we have that
\begin{equation}\label{10.13-eq1}
\begin{array}{ll}\ds
\int_0^L\|\sum_{m< |j|\leq
n}\big[\,\ov{\f_{j}'(x)}-c_j\phi_j'(x)\big]\f_{j}(y)\|^2dy\\
\ns\ds =\sum_{m< |j|\leq
n}\big|\ov{\f_{j}'(x)}-c_j\phi_j'(x)\big|^2\\
\ns \ds \leq 2\sum_{m< |j|\leq
n}\big|(1-c_j)\ov{\f_{j}'(x)}\big|^2 +
2\sum_{m< |j|\leq
n}|c_j|^2\big|\ov{\f_{j}'(x)}-\phi_j'(x)\big|^2.
\end{array}
\end{equation}
By means of \eqref{estimateLinftyfjx(x)}
and \eqref{11.7-eq22}, we find that
\begin{equation}\label{11.8-eq2}
\max_{x\in [0,L]}\sum_{m< |j|\leq
n}\big|(1-c_j)\ov{\f_{j}'(x)}\big|^2 \leq
C\sum_{m< |j|\leq n}\frac{\ln^2|j|}{j^2}.
\end{equation}
From \eqref{11.7-eq1} and
\eqref{8.30-eq15}, similarly to the proof
of \eqref{11.7-eq12}, we obtain that
\begin{equation}\label{11.8-eq3}
\begin{array}{ll}\ds
\q\max_{x\in[0,L]}\sum_{m< |j|\leq
n}|c_j|^2\big|\ov{\f_{j}'(x)}-\phi_j'(x)\big|^2\\
\ns\ds = \max_{x\in[0,L]}\sum_{m< |j|\leq
n}|c_j|^2\Big|\sum_{l\in\dbZ\setminus\{0,j\}}(i\mu_l-i\mu_{j}+\l)^{-1}
\frac{\l\f'_l(0)}{\f'_{j}(0)}\f_l'(x)\Big|^2\\
\ns\ds \leq C \sum_{m<|j|\leq
n}\sum_{l\in\dbZ\setminus\{0,j\}}\Big|\frac{l^2}{(j^3-l^3)j}\Big|^2
\\
\ns\ds \leq C\sum_{m<|j|\leq
n}\frac{\ln^2|j|}{j^2}.
\end{array}
\end{equation}
Combining \eqref{11.8-eq2} and
\eqref{11.8-eq3}, we get that
\begin{equation}\label{10.13-eq5}
\left\{ x\in (0,L) \mapsto
k^{(n)}_x(x,\cd)\in
L^2(0,L)\right\}_{n=1}^{+\infty} \mbox{ is
a Cauchy sequence in }C^0([0,L];L^2(0,L)),
\end{equation}
which shows that \eqref{regkx} holds.
Proceeding as in the proofs of
\eqref{kyL^2} and of \eqref{10.13-eq5}, one
gets that
\begin{equation}\label{cvkny}
\left\{ y\in (0,L) \mapsto
k^{(n)}_y(\cd,y)\in
L^2(0,L)\right\}_{n=1}^{+\infty} \mbox{ is
a Cauchy sequence in }C^0([0,L];L^2(0,L)),
\end{equation}
which also gives \eqref{regky}. Moreover
\eqref{11.2-eq2}, \eqref{11.8-eq4},
\eqref{defkn} and \eqref{cvkny} imply that
\begin{equation}\label{11.8-eq6}
k_y(\cd,0)=\lim_{n\to+\infty}k_y^{(n)}(\cd,0)=0
\;\mbox{ in } L^2(0,L).
\end{equation}
Similarly, one can show that
\begin{equation}\label{11.8-eq7}
k_y(\cd,L)=\lim_{n\to+\infty}k_y^{(n)}(\cd,L)=0
\;\mbox{ in } L^2(0,L).
\end{equation}
 From \eqref{11.8-eq6} and \eqref{11.8-eq7}, one has \eqref{ky0}.

Let us finally prove that $k(\cd,\cd)$
satisfies \eqref{deftransposition}.

First, it is clear that
\begin{equation}\label{11.11-eq1}
\begin{array}{ll}\ds
\q(\pa_{xxx}+\pa_x + \pa_{yyy}+\pa_y +
\l)(1-c_j)\ov{\f_j(x)}\f_j(y)\\
\ns\ds
=(1-c_j)(i\mu_j-i\mu_j+\l)\ov{\f_j(x)}\f_j(y)=(1-c_j)\l\ov{\f_j(x)}\f_j(y).
\end{array}
\end{equation}
From \eqref{defhatb} to \eqref{11.7-eq1},
one has
\begin{equation}\label{11.12-eq1}
\begin{array}{ll}\ds
\q(\pa_{xxx}+\pa_x + \pa_{yyy}+\pa_y +
\l)c_j[\ov{\f_j(x)}-\phi_j(x)]\f_j(y)\\
\ns\ds = -(\pa_{xxx}+\pa_x +
\pa_{yyy}+\pa_y +
\l)\sum_{k\in\dbZ\setminus\{0,j\}}c_j\frac{\l\f_k'(0)\ov{\f_k(x)}}{\f_j'(0)(i\mu_k-i\mu_j+\l)}\f_j(y)\\
\ns\ds =
-c_j\l\sum_{k\in\dbZ\setminus\{0,j\}}\frac{\f_k'(0)}{\f_j'(0)}\ov{\f_k(x)}\f_j(y)\\
\ns\ds =
c_j\l\ov{\f_j(x)}\f_j(y)-\frac{c_j\l}{\f_j'(0)}\sum_{k\in\dbZ\setminus\{0\}}
\f_k'(0)\ov{\f_k(x)}\f_j(y)\\
\ns\ds =
c_j\l\ov{\f_j(x)}\f_j(y)+\frac{c_j\l}{\f_j'(0)}\d'_{x=0}\otimes\f_j(y)\q
\mbox{ in }D(A)'\otimes L^2(0,L).
\end{array}
\end{equation}
From \eqref{defkn}, \eqref{11.11-eq1} and
\eqref{11.12-eq1}, we get in $D(A)'\otimes
L^2(0,L)$,
$$
\begin{array}{ll}\ds
k^{(n)}_{xxx}(x,y)+k^{(n)}_{x}(x,y)+k^{(n)}_{yyy}(x,y)+k^{(n)}_{y}(x,y)+\l
k^{(n)}(x,y)-\d'_{x=0}\otimes\sum_{0< |j|\leq n}\frac{c_j\l}{\f_j'(0)}\f_j(y)\\
\ns\ds= \l\sum_{0< |j|\leq
n}\ov{\f_j(x)}\f_j(y).
\end{array}
$$
Therefore, for any $\rho\in\cE\subset
D(A)\otimes L^2(0,L)$ (recall \eqref{def
cE} for the definition of $\cE$),  we have
\begin{equation}\label{7.30-eq1}
\begin{array}{ll}\ds
0\3n & \ds= \int_0^L \int_0^L \big[
k^{(n)}_{xxx}(x,y) + k^{(n)}_{x}(x,y)
+\!k^{(n)}_{yyy}(x,y) + k^{(n)}_{y}(x,y)
+ \!\l k^{(n)}(x,y)\\
\ns&\ds \q - \l\sum_{0< |j|\leq
n}\ov{\f_j(x)}\f_j(y)\big]\rho(x,y)dxdy -
\int_0^L\sum_{0< |j|\leq
n}\frac{\l}{\f_j'(0)}\d'_{x=0}\big(\rho(x,y)\big)\f_j(y)dy\\
\ns&\ds = -\int_0^L \int_0^L  \big[
\rho_{xxx}(x,y) + \rho_{x}(x,y) +
\rho_{yyy}(x,y) + \rho_{y}(x,y) +
\l\rho(x,y)
\big]k^{(n)}(x,y)dxdy\\
\ns&\ds \q-\int_0^L
k_y^{(n)}(x,L)\rho(x,L)dx + \int_0^L
k_y^{(n)}(x,0)\rho(x,0)dx
\\
\ns&\ds \q - \l\int_0^L \int_0^L
\rho(x,y)\sum_{0< |j|\leq
n}\ov{\f_j(x)}\f_j(y) dxdy.
\end{array}
\end{equation}
By \eqref{11.8-eq6}, \eqref{11.8-eq7} and
letting $n\to+\infty$ in \eqref{7.30-eq1},
we obtain that
\begin{equation}\label{7.30-eq1-new}
\begin{array}{ll}\ds
\int_0^L \int_0^L  \big[ \rho_{xxx}(x,y) +
\rho_{x}(x,y) + \rho_{yyy}(x,y) +
\rho_{y}(x,y) + \l\rho(x,y)
\big]k(x,y)dxdy\\
\ns \ds  + \int_0^L \rho(y,y) dy =0.
\end{array}
\end{equation}
 This concludes Step 2 and
therefore the proof of Lemma
\ref{well-posed lm}
\begin{remark}
\label{remarkrealvalued} Note that
$k(\cd,\cd)$ is a real valued function
since $\bar k(\cd,\cd)$ is also a solution
of \eqref{system3} in $\cG$ and therefore,
by the uniqueness proved in Step 1, we must
have $\bar k(\cd,\cd)=k(\cd,\cd)$.
\end{remark}

\vspace{0.3cm}

\section{Invertibility of $I-K$}
\label{secinvert} We define  a bounded
linear operator $K:L^2(0,L)\to L^2(0,L)$ by
$$
K(v)(x)\=\int_0^L k(x,y)v(y)dy\q\mbox{ for
every }v\in L^2(0,L).
$$
Note that, by Remark
\ref{remarkrealvalued}, if $v$ is a real
valued function, then $Kv$ is also a real
valued function. The goal of this section
is to prove the following lemma.
\begin{lemma}\label{invertible lm}
$I-K$ is an invertible operator.
\end{lemma}

{\it Proof of Lemma~\ref{invertible lm}}\,:
Since $ k(\cd,\cd)\in L^2((0,L)\times
(0,L))$, we get that $K$ is a compact
operator. Further, since $k\in
H^{1}_0((0,L)\times (0,L)) $, we know that
$K$ is a continuous linear map from
$L^2(0,L)$ into $H^{1}_0(0,L)$.  Denote by
$K^*$ the adjoint operator of $K$. Then, it
is easy to see that
$$
K^* (v)(x) = \int_0^L k^*(x,y)v(y)dy
\q\mbox{ for any } v\in L^2(0,L),
$$
where $k^*$ is defined by
\begin{equation}\label{defk*}
  k^*(x,y)\=k(y,x),\, (x,y)\in (0,L)\times(0,L).
\end{equation}
 From \eqref{system3} and \eqref{defk*}, one gets
\begin{equation}\label{system3.31}
\left\{
\begin{array}{ll}\ds
k^*_y  +  k^*_{yyy}  +  k^*_x +
 k^*_{xxx} +\l  k^* =\l\d(y-x) &\mbox{ in
}
(0,L)\times (0,L),\\
\ns\ds  k^*(0,y) =  k^*(L,y) =0 &\mbox{ on } (0,L),\\
\ns\ds  k^*_x(0,y) =  k^*_x(L,y) =0 &\mbox{ on } (0,L),\\
\ns\ds  k^*(x,0) =  k^*(x,L)=0 &\mbox{ on }
(0,L).
\end{array}
\right.
\end{equation}
Further, from    \eqref{defk*} again, we
have the following regularity for
$k^*(\cd,\cd)$:
\begin{gather}
\label{9.5-eq5-1*} y\in (0,L)\mapsto
k_y^*(\cd,y)\in L^2(0,L) \text{ is in }
C^0([0,L]; L^2(0,L)),
\\
\label{9.5-eq5-2*} x\in (0,L)\mapsto
k_x^*(x,\cdot)\in L^2(0,L) \text{ is in }
C^0([0,L]; L^2(0,L)).
\end{gather}
Let us point out that, from the
\eqref{system3.31} and \eqref{9.5-eq5-2*},
we know

\vspace{-0.35cm}

\begin{equation}\label{9.5-eq2}
v\in K^*(L^2(0,L))\Rightarrow \left(v\in
C^1([0,L]),\;
v(0)=v(L)=v_x(0)=v_x(L)=0\right).
\end{equation}

We claim that the spectral radius $r(K^*)$
of $K^*$ equals $0$. Otherwise, since $K^*:
L^2(0,L)\to L^2(0,L)$ is as $K$ a compact
linear operator, there is a nonzero
eigenvalue $\a$ of $K^*$. Then, there
exists a positive integer $n_0$ such that
\begin{equation}
\label{n0+1} \Ker(K^* - \a
I)^{n_0+1}=\Ker(K^* - \a I)^{n_0}.
\end{equation}
Let
$$
\cF\=\Ker(K^* - \a I)^{n_0}.
$$
It is a finite dimensional space.

Since $\a\not =0$, $\cF\subset
K^*(L^2(0,L))$, and, with \eqref{9.5-eq2},
\begin{equation}
\label{vsurbord} \cF\subset C^1([0,L])
\text{ and } v(0)=v_x(0)=v(L)=v_x(L)=0, \,
\forall\, v \in \cF.
\end{equation}
Note that, by   the fact that $k^*\in
H^1_0((0,L)\times(0,L))$,
\begin{equation}
\label{k*H-1} \text{$K^*$ can be extended
to be a continuous linear map from
$H^{-1}(0,L)$ into $L^2(0,L)$.}
\end{equation}
Let us denote by $\wt K^*$ this extension
and remark that, if $u\in H^{-1}(0,L)$ is
such that $\wt K^* u =\a u$, then $u\in
L^2(0,L)$. Thus, we see that $\Ker(\wt
K^*-\a I)\subset L^2(0,L)$ and $\Ker(\wt
K^*-\a I)=\Ker(K^*-\a I)$. Similarly, we
have
$$
\left\{
\begin{array}{ll}
\ds\Ker(\wt K^*-\a I)^{n_0}\subset
L^2(0,L),\\
\ns\ds  \Ker(\wt K^*-\a I)^{n_0}=\cF, \\
\ns\ds\Ker(\wt K^*-\a I)^{n_0}=\Ker(\wt
K^*-\a I)^{n_0+1}.
\end{array}
\right.
$$
By \eqref{system3.31},
\begin{gather}
\label{commutation} K^* (\pa_{xxx}+\pa_x) v
= (\pa_{xxx}+\pa_x) K^* v -\l K^* v + \l v,
\, \forall\, v\in C^\infty_0(0,L).
\end{gather}
From \eqref{9.5-eq2} and
\eqref{commutation} we also get that $ K^*$
can be extended to be a continuous linear
map from $H^{-3}(0,L)$ into $H^{-2}(0,L)$.
This, together with \eqref{k*H-1} and an
interpolation argument, shows that $ K^*$
can be extended to be a continuous linear
map from $H^{-2}(0,L)$ into $H^{-1}(0,L)$.
We denote by $\wh K^*$ this extension.
Then, as for $ \wt K^*$,  one has
\begin{equation}\label{9.5-eq4}
\left\{
\begin{array}{ll}
\ds\Ker(\wh K^*-\a I)^{n_0} =\Ker(K^*-\a
I)\in
L^2(0,L),\\
\ns\ds  \Ker(\wh K^*-\a I)^{n_0}=\cF, \\
\ns\ds\Ker(\wh K^*-\a I)^{n_0}=\Ker(\wh
K^*-\a I)^{n_0+1}.
\end{array}
\right.
\end{equation}
Using \eqref{commutation}, a density
argument and \eqref{k*H-1} we get that
\begin{equation}\label{commutationwithA}
\wh K^* (\pa_{xxx}+\pa_x) v =
(\pa_{xxx}+\pa_x) \wh K^* v -\l \wh K^* v +
\l v, \, \forall v\in \cF.
\end{equation}

From \eqref{commutationwithA},
\eqref{vsurbord}, and induction on $n$, one
gets that
\begin{equation}\label{K*n}
(\wh K^*)^n (\pa_{xxx}+\pa_x) v=
(\pa_{xxx}+\pa_x) (\wh K^*)^n v - n\lambda
(\wh K^*)^n v + n\l(\wh K^*)^{n-1} v, \,
\forall v\in \cF.
\end{equation}
and therefore, for every polynomial $P$,
\begin{equation}\label{K*poly}
P(\wh K^*) (\pa_{xxx}+\pa_x) v =
(\pa_{xxx}+\pa_x) P(\wh K^*)v - \lambda
P'(\wh K^*)\wh K^* v + \lambda P'(\wh K^*)
v, \, \forall v\in \cF.
\end{equation}
By virtue of \eqref{9.5-eq4},
\eqref{vsurbord} and \eqref{K*poly} with
$P(X)\=(X-\a)^{n_0+1}$, we see that
$(\pa_{xxx}+\pa_x)\cF\subset \cF$. Since
$\cF$ is finite dimensional, this implies
that $(\pa_{xxx}+\pa_x)$ has an
eigenfunction in $\cF$, that is, there
exist  $\mu\in \dbC$ and
$\xi\in\cF\setminus\{0\}$ such that

$$
\left\{
\begin{array}{ll}\ds
(\pa_{xxx}+\pa_x)\xi = \mu\xi &\mbox{ in } (0,L),\\
\ns\ds \xi(0)=\xi(L)=\xi_x(0)=\xi_x(L)=0.
\end{array}
\right.
$$
But, by \cite[Lemma 3.5]{Rosier1} again,
this is impossible since $L\notin\cN$.
Then, we know $r(K^*)=0$, which implies
that the spectral radius $r(K)$ of $K$ is
zero. Hence, we know that the real number
$1$ belongs to the resolvent set of $K$,
which completes the proof of
Lemma~\ref{invertible lm}.


\section{Proof of Theorem \ref{th1}}
\label{secproofth}


This section is addressed to a proof of
Theorem \ref{th1}.

{\it Proof of Theorem \ref{th1}}\,: Let
$T>0$, which will be given later. Consider
the following equation

\begin{equation}\label{12.22-eq6}
\left\{
\begin{array}{ll}\ds
v_{1,t} + v_{1,x} + v_{1,xxx} + v_1v_{1,x}
= 0 &\mbox{ in }
[0,T]\times (0,L),\\
\ns\ds v_1(t,0) = v_1(t,L)=0 &\mbox{ on }
[0,T],\\
\ns\ds v_{1,x}(t,L)=\int_0^L
k_x(L,y)v_1(t,y)dy&\mbox{ on } [0,T],
\\
\ns\ds v_1(0)=v^0 & \mbox{ in } (0,L).
\end{array}
\right.
\end{equation}
By Theorem \ref{well}, we know that there
is an $r_T>0$ such that for all $v^0\in
L^2(0,L)$ with $|v^0|_{L^2(0,L)}\leq r_T$,
the equation \eqref{12.22-eq6} admits a
unique solution $v_1\in
C^0([0,T];L^2(0,L))\cap
L^2(0,T;H^1_0(0,L))$. Moreover, there is a
constant $C_T>0$ such that
\begin{equation}\label{12.22-eq2}
|v_1|_{C^0([0,T];L^2(0,L))\cap
L^2(0,T;H^1_0(0,L))} \leq
C_T|v^0|_{L^2(0,L)}.
\end{equation}

Let $w_1 = (I-K)v_1$. Then, we have that
\begin{equation}\label{eq1}
\begin{array}{ll}\ds
w_{1,t}(t,x)\3n &\ds =v_{1,t}(t,x) -
\int_0^L
k(x,y)v_{1,t}(t,y)dy\\
\ns&\ds = v_{1,t}(t,x) + \int_0^L
k(x,y)\big[v_{1,y}(t,y) + v_{1,yyy}(t,y)+v_1(t,y)v_{1,y}(t,y)\big]dy\\
\ns&\ds = v_{1,t}(t,x) +  \big[k(x,y)
v_1(t,y)\big]\big|_{y=0}^{y=L} - \int_0^L
k_y(x,y) v_1(t,y)dy \\
\ns&\ds \q +  \big[ k(x,y)v_{1,yy}(t,y)
\big]\big|_{y=0}^{y=L}\! - \! \big[
k_y(x,y)v_{1,y}(t,y)
\big]\big|_{y=0}^{y=L}\! + \!\big[
k_{yy}(x,y)v_1(t,y)
\big]\big|_{y=0}^{y=L} \\
\ns&\ds \q- \!\int_0^L\!\!
k_{yyy}(x,y)v_1(t,y)dy \!+\!
\frac{1}{2}\big[ k(x,y)v_1(t,y)^2
\big]\big|_{y=0}^{y=L}\! -\!
\frac{1}{2}\!\int_0^L\!\!
k_y(x,y)v_1(t,y)^2dy
 \\
\ns&\ds = v_{1,t}(t,x) - \int_0^L
\big[k_y(x,y)+k_{yyy}(x,y)\big] v_1(t,y)dy
- \frac{1}{2}\int_0^L k_y(x,y)v_1(t,y)^2dy,
\end{array}
\end{equation}
\begin{equation}\label{eq2}
w_{1,x}(t,x) = v_{1,x}(t,x) - \int_0^L
k_x(x,y)v_1(t,y)dy
\end{equation}
and
\begin{equation}\label{eq3}
w_{1,xxx}(t,x) = v_{1,xxx}(t,x) - \int_0^L
k_{xxx}(x,y)v_1(t,y)dy.
\end{equation}
Therefore, for a given $\l\in\dbR$,
utilizing the fact the $v_1$ solves
\eqref{12.22-eq6}, we obtain that
\begin{equation}\label{eq4}
\begin{array}{ll}\ds
\q \vspace{-0.2cm} w_{1,t}(t,x) +
w_{1,x}(t,x) + w_{1,xxx}(t,x) + \l
w_1(t,x)+  \frac{1}{2}\int_0^L k_y(x,y)v_1(t,y)^2dy \\
\ns\ds \vspace{-0.2cm} = \l
v_1(t,x)-\int_0^L v_1(t,y)\big[ k_y(x,y) +
k_{yyy}(x,y) + k_x(x,y) + k_{xxx}(x,y) + \l
k(x,y) \big] dy\\
\ns\ds = -\int_0^L v_1(t,y)\big[ k_y(x,y) +
k_{yyy}(x,y) + k_x(x,y) + k_{xxx}(x,y) + \l
k(x,y) -\l \d(x-y) \big] dy \\
\ns\ds = -v_1(t,x)v_{1,x}(t,x).
\end{array}
\end{equation}
Hence, if we   take the feedback control
$F$ as
\begin{equation}\label{deffeedback-varphi}
F(\varphi)\=\int_0^L k_x(L,y)\varphi(y)dy,
\q \forall \varphi \in L^2(0,L),
\end{equation}
then we get that $w_1$ solves
\begin{equation}\label{6.9-eq12}
\left\{ \begin{array}{ll}\ds w_{1,t}+
w_{1,x} + w_{1,xxx} + w_1w_{1,x} + \l w_1 \\
\ns\ds   = - v_1v_{1,x} - \frac{1}{2}
\int_0^L k_y(x,y)v_1(t,y)^2dy &\mbox{in }
[0,T]\times (0,L),\\
\ns\ds  w_1(t,0)=w_1(t,L)=0 &\mbox{on }[0,T],\\
\ns\ds  w_{1,x}(t,L)=0 &\mbox{on } [0,T].
\end{array}
\right.
\end{equation}
Thus, we see that
\begin{equation}\label{6.9-eq13}
\begin{array}{ll}\ds
\q\frac{d}{dt}\int_0^L| w_1(t,x)|^2 dx \\
\ns\ds = -|w_{1,x}(t,0)|^2 -
2\l\int_0^L|w_1(t,x)|^2 dx -
\int_0^L w_1(t,x)v_1(t,x)v_{1,x}(t,x)dx\\
\ns\ds \q - \frac{1}{2}\int_0^L  w_1(t,x)
\[\int_0^L k_y(x,y)v_1(t,y)^2dy\]dx.
\end{array}
\end{equation}
Now we estimate the third and  fourth terms
in the right hand side of \eqref{6.9-eq13}.
First,
\begin{equation}\label{12.22-eq3}
\begin{array}{ll}\ds
\q\Big|\int_0^L
w_1(t,x)v_1(t,x)v_{1,x}(t,x)dx\Big|
\\
\ns\ds = \Big|\int_0^L
\[v_1(t,x)-\int_0^L k(x,y)v_1(t,y)dy\]v_1(t,x)v_{1,x}(t,x)dx \Big|
\\
\ns\ds = \Big|\int_0^L \(\int_0^L
k(x,y)v_1(t,y)dy\) v_1(t,x)v_{1,x}(t,x)dx
\Big|
\\
\ns\ds = \frac{1}{2}\Big|\int_0^L
\(\int_0^L k_x(x,y)v_1(t,y)dy\)
v_1(t,x)^2dx \Big|
\\
\ns\ds \leq
\frac{1}{2}\Big|\int_0^L|k_x(\cd,y)|^2dy\Big|_{L^\infty(0,L)}\(\int_0^L
v_1(t,x)^2dx\)^{\frac{3}{2}}
\\
\ns\ds \leq |(I-K)^{-1}|^3_{\cL(L^2(0,L))}
\frac{1}{2}\Big|\int_0^L|k_x(\cd,y)|^2dy\Big|_{L^\infty(0,L)}\(\int_0^L
w^2_1(t,x)dx\)^{\frac{3}{2}}.
\end{array}
\end{equation}
Next,
\begin{equation}\label{6.9-eq14}
\!\!\!\begin{array}{ll}\ds \q \Big|\int_0^L
w_1(t,x)\int_0^L
k_y(x,y)v_1(t,y)^2dy dx  \Big| \\
\ns\ds = \Big| \int_0^L \[
v_1(t,x)-\int_0^L k(x,y)v_1(t,y)dy
\]\(\int_0^L
k_y(x,y)v_1(t,y)^2dy\)dx\Big|\\
\ns \ds \leq \!\Big|\!\int_0^L \!\!v_1(t,x)
\(\!\int_0^L\!\! k_y(x,y)v_1(t,y)^2dy
\)dx\Big| \!+\! \Big|\! \int_0^L\!\!
\(\!\int_0^L\!\! k(x,y)v_1(t,y)dy
\int_0^L\!\! k_y(x,y)v_1(t,y)^2dy
\)dx\Big|.
\end{array}
\end{equation}
The first term in the right hand side of
\eqref{6.9-eq14} satisfies that
\begin{equation}\label{10.15-eq1}
\begin{array}{ll}\ds
\q\Big| \int_0^L v_1(t,x) \(\int_0^L
k_y(x,y)v_1(t,y)^2dy\)dx\Big|\\
\ns\ds = \Big| \int_0^L v_1(t,y)^2
\(\int_0^L
k_y(x,y)v_1(t,x) dx\) dy \Big|\\
\ns\ds \leq \int_0^L v_1(t,y)^2 \(\int_0^L
|k_y(x,y)|^2 dx\)^{\frac{1}{2}}dy
\(\int_0^L
v_1(t,x)^2dx\)^{\frac{1}{2}}\\
\ns\ds \leq  \Big|\(\int_0^L |k_y(x,\cd)|^2
dx\)^{\frac{1}{2}}\Big|_{L^\infty(0,L)}
\(\int_0^L v_1(t,x)^2dx\)^{\frac{3}{2}}\\
\ns\ds \leq
|(I-K)^{-1}|_{\cL(L^2(0,L))}^{3}\Big|\(\int_0^L
|k_y(x,\cd)|^2
dx\)^{\frac{1}{2}}\Big|_{L^\infty(0,L)}
\(\int_0^L w_1(t,x)^2dx\)^{\frac{3}{2}}.
\end{array}
\end{equation}
The second term in the right hand side of
\eqref{6.9-eq14} satisfies that
\begin{equation}\label{10.15-eq2}
\begin{array}{ll}\ds
\q\Big|\int_0^L \(\int_0^L k(x,z)v_1(t,z)dz
\int_0^L
k_y(x,y)v_1(t,y)^2dy\)dx\Big|\\
\ns\ds = \Big| \int_0^L v_1(t,y)^2 \int_0^L
k_y(x,y)\(\int_0^L k(x,z)v_1(t,z)dz\)dx dy \Big|\\
\ns\ds \leq \int_0^L v_1(t,y)^2 \(\int_0^L
|k_y(x,y)|^2 dx\)^{\frac{1}{2}}dy
\[\int_0^L
\(\int_0^L k(x,z)v_1(t,z)dz\)^2dx\]^{\frac{1}{2}}\\
\ns\ds \leq  \Big|\(\int_0^L |k_y(x,\cd)|^2
dx\)^{\frac{1}{2}}\Big|_{L^\infty(0,L)}\(\int_0^L
|k(x,y)|^2 dxdy\)^{\frac{1}{2}}
\(\int_0^L v_1(t,x)^2dx\)^{\frac{3}{2}}\\
\ns\ds \leq
|(I\!-\!K)^{-1}|_{\cL(L^2(0,L))}^{3}\Big|\(\!\int_0^L\!
\!|k_y(x,\cd)|^2
dx\)^{\frac{1}{2}}\Big|_{L^\infty(0,L)}\!\(\!\int_0^L\!\!\!
|k(x,y)|^2 dxdy\)^{\frac{1}{2}}
\(\!\int_0^L\!\!\!
w_1(t,x)^2dx\)^{\frac{3}{2}}.
\end{array}
\end{equation}

Let
$$
\begin{array}{ll}\ds
\wh C \3n&\ds=
\frac{1}{2}|(I-K)^{-1}|^3_{\cL(L^2(0,L))}
\Big|\int_0^L|k_x(\cd,y)|^2dy\Big|_{L^\infty(0,L)}
\\
\ns&\ds \q +
|(I-K)^{-1}|_{\cL(L^2(0,L))}^{3}\Big|\(\int_0^L
|k_y(x,\cd)|^2
dx\)^{\frac{1}{2}}\Big|_{L^\infty(0,L)}
\\
\ns&\ds \q +
|(I-K)^{-1}|_{\cL(L^2(0,L))}^{3}\Big|\(\int_0^L
|k_y(x,\cd)|^2
dx\)^{\frac{1}{2}}\Big|_{L^\infty(0,L)}\(\int_0^L
|k(x,y)|^2 dxdy\)^{\frac{1}{2}}.
\end{array}
$$
From  \eqref{6.9-eq13}--\eqref{10.15-eq2},
we get that
\begin{equation}\label{12.22-eq4}
\frac{d}{dt}\int_0^L | w_1(t,x)|^2 dx  \leq
 - 2\l \int_0^L
|w_1(t,x)|^2 dx +\wh C\( \int_0^L
w_1(t,x)^2dx\)^{\frac{3}{2}}.
\end{equation}
For a given $\l>0$, we know that there is a
$\d_1>0$ such that, if $|w_1(0)|_{L^2(0,L)}
\leq \d_1$, then
\begin{equation*}
\wh C\( \int_0^L
w_1(t,x)^2dx\)^{\frac{1}{2}} \leq  \l,\quad
\forall\, t\in [0,T].
\end{equation*}
This, together with \eqref{12.22-eq4},
implies that
\begin{equation*}
\frac{d}{dt}\int_0^L |w_1(t,x)|^2 dx \leq
- \l\int_0^L|w_1(t,x)|^2 dx,\quad \forall\,
t\in [0,T].
\end{equation*}
In particular
\begin{equation}\label{6.9-eq15}
|w_1(t,\cd)|_{L^2(0,L)} \leq e^{-\frac{\l
t}{2}}|w_1(0,\cd)|_{L^2(0,L)},
\q\forall\,t\in [0,T].
\end{equation}
Then, from Lemma \ref{invertible lm} and
\eqref{6.9-eq15}, we get that  if
$$
|v_1(0,\cd)|_{L^2(0,L)}\leq
\min\{|I-K|^{-1}_{\cL(L^2(0,L))}\d_1,
r_T\},
$$
then
\begin{equation}\label{6.9-eq16}
|v_1(t,\cd)|_{L^2(0,L)} \leq
|(I-K)^{-1}|_{\cL(L^2(0,L))}|I-K|_{\cL(L^2(0,L))}e^{-\frac{\l
t}{2}}|v_1(0,\cd)|_{L^2(0,L)},\q \forall\,
t\in [0,T].
\end{equation}
Now we choose $T>0$ such that
$$
e^{-\frac{\l T}{2}} \leq
|(I-K)^{-1}|_{\cL(L^2(0,L))}|I-K|_{\cL(L^2(0,L))}.
$$
From \eqref{6.9-eq16}, we find that
$|v_1(T)|_{L^2(0,L)}\leq
|v_1(0)|_{L^2(0,L)} \leq r_T$. Thus, by
Theorem \ref{well}, we know that the
following equation is well-posed.
\begin{equation}\label{12.22-eq5}
\left\{
\begin{array}{ll}\ds
v_{2,t} + v_{2,x} + v_{2,xxx} + v_2v_{2,x}
= 0 &\mbox{ in }
[0,T]\times (0,L),\\
\ns\ds v_2(t,0) = v_2(t,L)=0 &\mbox{ on }
[0,T],\\
\ns\ds v_{2,x}(t,L)=\int_0^L
k_x(L,y)v_2(t,y)dy&\mbox{ on } [0,T], \\
\ns\ds v_2(0)= v_1(T) &\mbox{ in } (0,L).
\end{array}
\right.
\end{equation}
Further, by a similar argument, we can find
that $w_2\=(I-K)v_2$ satisfies
$$
|w_2(t,\cd)|_{L^2(0,L)} \leq e^{-\frac{\l
t}{2}}|w_2(0,\cd)|_{L^2(0,L)}, \q\forall\,
t\in [0,T]
$$
and
$$
|v_2(T,\cd)|_{L^2(0,L)}\leq
|v_2(0,\cd)|_{L^2(0,L)}.
$$
Then, we can define $v_3$ and $w_3$ in a
similar manner. By induction, we can find
$v_n\in C^0([0,T];L^2(0,L))\cap
L^2(0,T;H^1_0(0,L))$ ($n>1$), which solves
\begin{equation}\label{12.22-eq6.1}
\left\{
\begin{array}{ll}\ds
v_{n,t} + v_{n,x} + v_{n,xxx} + v_nv_{n,x}
= 0 &\mbox{ in }
[0,T]\times (0,L),\\
\ns\ds v_n(t,0) = v_n(t,L)=0 &\mbox{ on }
[0,T],\\
\ns\ds v_{n,x}(t,L)=\int_0^L
k_x(L,y)v_n(t,y)dy&\mbox{ on } [0,T], \\
\ns\ds v_n(0)= v_{n-1}(T) &\mbox{ in }
(0,L).
\end{array}
\right.
\end{equation}
Further, we have that $w_n = (I-K)v_n$
satisfies
\begin{equation}\label{12.22-eq7}
|w_n(t,\cd)|_{L^2(0,L)} \leq e^{-\frac{\l
t}{2}}|w_{n}(0,\cd)|_{L^2(0,L)}=
e^{-\frac{\l
t}{2}}|w_{n-1}(T,\cd)|_{L^2(0,L)}
\end{equation}
and
$$
|v_n(T,\cd)|_{L^2(0,L)} \leq
|v_{n}(0,\cd)|_{L^2(0,L)}=|v_{n-1}(T,\cd)|_{L^2(0,L)}.
$$
Now we put
$$
v(t+(n-1)T,x) = v_n(t,x),\q
w(t+(n-1)T,x)=w_n(t,x) \q\mbox{ for }
(t,x)\in [0,T)\times [0,L].
$$
Then, it is an easy matter to see that $v$
is a solution to \eqref{csystem1} and
$w=(I-K)v$. From \eqref{12.22-eq7}, we get
that
\begin{equation}\label{12.22-eq8}
|w(t,\cd)|_{L^2(0,L)}\leq e^{-\frac{\l
t}{2}}|w(0,\cd)|_{L^2(0,L)},\q\forall\,
t\geq 0.
\end{equation}
This, together with $w=(I-K)v$, implies
that for all $t\geq 0$,
$$
|v(t,\cd)|_{L^2(0,L)}\leq e^{-\frac{\l
t}{2}}|(I-K)^{-1}|_{\cL(L^2(0,L))}|I-K|_{\cL(L^2(0,L))}|v(0,\cd)|_{L^2(0,L)}\leq
Ce^{-\frac{\l t}{2}}|v(0,\cd)|_{L^2(0,L)}.
$$
Let
$\d_0=\min\{|(I-K)|^{-1}_{\cL(L^2(0,L))}\d_1,
r_T\}$. Then, we know that for any $v^0\in
L^2(0,L)$ with $|v^0|_{L^2(0,L)}\leq \d_0$,
the equation \eqref{csystem1} with the
initial datum $v^0$ admits a solution
$$
v \in C^0([0,+\infty);L^2(0,L))\cap
L^2_{loc}(0,+\infty;H^1_0(0,L)).
$$
Further, we have that
\begin{equation}\label{12.22-eq9}
|v(t,\cd)|_{L^2(0,L)}\leq Ce^{-\frac{\l
t}{2}}|v(0,\cd)|_{L^2(0,L)},\q\forall\,
t\geq 0.
\end{equation}
%


\appendix
\section{Appendix}

This section is devoted to a proof of
Theorem \ref{well}. Before giving the
proof, we first recall the following useful
results.

Let $T>0$. We consider the following
linearized KdV equation with
non-homogeneous boundary condition.
\begin{equation}\label{12.19-eq1}
\left\{
\begin{array}{ll}\ds
u_t + u_{xxx} + u_x =\tilde h &\mbox{ in }
[0,T]\times (0,L),\\
\ns\ds u(t,0)=u(t,L)=0 &\mbox{ on } [0,T],\\
\ns\ds u_x(t,L) = h(t) &\mbox{ on } [0,T],\\
\ns\ds u(0)=u^0 &\mbox{ in } (0,L).
\end{array}
\right.
\end{equation}
Here $u^0\in L^2(0,L)$,  $h\in L^2(0,T)$,
$\tilde h\in L^1(0,T;L^2(0,L))$. In
\cite{Rosier1} (see also
\cite{2014-Cerpa-MCRF}), the author proved
the following results:

\begin{lemma}\label{lm3}
Let $u^0\in L^2(0,L)$. There exists a
unique solution  $u\in
C^0([0,T];L^2(0,L))\cap
L^2(0,T;H^1_0(0,L))$ of \eqref{12.19-eq1}
such that $u(0,\cdot)=u^0(\cdot)$.
Moreover, there is a $C_1>0$, independent
of $h\in L^2(0,T)$, $\tilde h\in
L^1(0,T;L^2(0,L))$ and $u^0\in L^2(0,L)$,
such that
\begin{equation} \label{est1}
|u|^2_{C^0([0,T];L^2(0,L))} +
|u|^2_{L^2(0,T;H^1_0(0,L))} \leq
C_1\big(|u^0|_{L^2(0,L)}^2 +
|h|_{L^2(0,T)}^2 + |\tilde
h|_{L^1(0,T;L^2(0,L))}^2\big).
\end{equation}
\end{lemma}
\begin{lemma}\label{lm4}
Let $z\in L^2(0,T;H^1(0,L))$. Then $zz_x\in
L^1(0,T;L^2(0,L))$ and the map $z\in
L^2(0,T;H^1(0,L))\mapsto zz_x\in
L^1(0,T;L^2(0,L))$ is continuous.
\end{lemma}

In \cite{CorCr}, the following result is
proved.

\begin{lemma}\label{lm5}\cite[Proposition 15]{CorCr}
 There exists $C>0$ such
that for every $\tilde u_0,\hat u_0\in
L^2(0,L)$ and  $\tilde h,\hat h \in
L^2(0,T)$ for which there exist solution
$\tilde u$ and $\hat u$ in
$C^0([0,T];L^2(0,L))\cap
L^2(0,T;H^1_0(0,L))$ of
\begin{equation}\label{eqtildeu-A}
\left\{
\begin{array}{ll}\ds
\tilde u_t + \tilde u_x + \tilde u_{xxx} +
\tilde u\tilde u_x =0 &\mbox{ in }
[0,T]\times (0,L),\\
\ns\ds \tilde u(t,0)=\tilde u(t,L)=0 &\mbox{ on } [0,T],\\
\ns\ds \tilde u_x(t,L) = \tilde h(t) &\mbox{ on } [0,T],\\
\ns\ds \tilde u(0,x) = \tilde u^0(x)
&\mbox{ in } (0,L),
\end{array}
\right.
\end{equation}
and of
\begin{equation}\label{eqhatu-A}
\left\{
\begin{array}{ll}\ds
\hat u_t + \hat u_x + \hat u_{xxx} + \hat
u\hat u_x =0 &\mbox{ in }
[0,T]\times (0,L),\\
\ns\ds \hat u(t,0)=\hat u(t,L)=0 &\mbox{ on } [0,T],\\
\ns\ds \hat u_x(t,L) = \hat h(t) &\mbox{ on } [0,T],\\
\ns\ds \hat u(0,x) = \hat u^0(x) &\mbox{ in
} (0,L),
\end{array}
\right.
\end{equation}
one has the following inequality
\begin{equation}\label{1.4-eq2}
\begin{array}{ll}\ds
\q\int_0^L|\tilde u(t,x)-\hat u(t,x)|^2dx
\\
\ns\ds \leq e^{C\big(1+|\tilde
u|_{L^2(0,T;H^1(0,L))}^2 + |\hat
u|_{L^2(0,T;H^1(0,L))}\big)^2}\[\int_0^L|\tilde
u^0(x) - \hat u^0(x)|^2 dx +
\int_0^t|\tilde h(s) - \hat h(s)|^2 ds\]
\end{array}
\end{equation}
for all $t\in [0,T]$.
\end{lemma}

\begin{remark}
In \cite{CorCr}, the term $\int_0^t|\tilde
h(s) - \hat h(s)|^2 ds$ in the right hand
side of \eqref{1.4-eq2} is $\int_0^T|\tilde
h(s) - \hat h(s)|^2 ds$. However, one can
see that the proof of \cite{CorCr} also
gives \eqref{1.4-eq2}.
\end{remark}

{\it Proof of Theorem \ref{well}}\,:

{\bf Uniqueness of the solution}

Assume that $v_1$ and $v_2$ are two
solutions of \eqref{csystem1}, by Lemma
\ref{lm5}, we know that for any $t\in
[0,T]$,
\begin{equation}\label{1.4-eq3}
\begin{array}{ll}\ds
\q\int_0^L |v_1(t,x)-v_2(t,x)|^2 dx
\\
\ns\ds \leq e^{C\big( 1+
|v_1|_{L^2(0,T;H^1(0,L))}^2 +
|v_2|_{L^2(0,T;H^1(0,L))}^2
\big)}\int_0^t|F(v_1)-F(v_2)|^2 ds
\\
\ns\ds \leq Ce^{C\big( 1+
|v_1|_{L^2(0,T;H^1(0,L))}^2 +
|v_2|_{L^2(0,T;H^1(0,L))}^2
\big)}\int_0^t\int_0^L
|v_1(s,x)-v_2(s,x)|^2dxds.
\end{array}
\end{equation}
This, together with the Gronwall
inequality, implies that $v_1=v_2$ in
$[0,T]\times (0,L)$.

\vspace{0.4cm}

{\bf Existence of the solution}

 Let us extend
$\tilde h$ and $h$ to be a function on
$(0,+\infty)\times (0,L)$ and $(0,+\infty)$
by setting them to be equal to zero on
$(T,+\infty)\times (0,L)$ and
$(T,+\infty)$, respectively. Denote by
$||F||$ the norm of the continuous linear
map $F:L^2(0,L)\to \dbR$. Set
\begin{gather}
\label{defT1}
T_1\=\min\Big\{\frac{1}{2C_1||F||^{2}},T\Big\}.
\end{gather}
Let $u\in C^0([0,T_1];L^2(0,L))$. We know
that
$$
h(\cd)\=  F(u(\cd)) \in L^2(0,T_1).
$$
Hence, for $v^0$ given in $L^2(0,L)$, we
can define a map
$$\cJ: C^0([0,T_1];L^2(0,L))\to C^0([0,T_1];L^2(0,L))$$
as follows:
\noindent $\cJ(u)=v$, where $v\in
C^0([0,T_1];L^2(0,L))\cap
L^2(0,T;H^1_0(0,L))$ solves
\eqref{12.19-eq1} with  $h(\cd) =
F(u(\cd))$ and  $v(0,\cdot)=v^0(\cdot)$.

 For $\tilde u, \hat u\in
C^0([0,T_1];L^2(0,L))$, one has, using in
particular \eqref{est1} and \eqref{defT1},
$$
\begin{array}{ll}\ds
\q|\cJ(\tilde u)-\cJ(\hat
u)|_{C^0([0,T_1];L^2(0,L))}
\\
\ns\ds \leq
\sqrt{C_1}\(\int_0^{T_1}\big|F(\tilde
u(t,\cd)) - F(\hat u(t,\cd)) \big|^2
dt\)^{\frac{1}{2}}
\\
\ns\ds \leq
\sqrt{C_1}||F||\(\int_0^{T_1}\int_0^L
|\tilde u(t,y)- \hat u(t,y)|^2
dydt\)^{\frac{1}{2}}
\\
\ns\ds \leq \sqrt{T_1}\sqrt{C_1} ||F||
|\tilde u-
\hat u|_{C^0([0,T_1];L^2(0,L))} \\
\ns\ds \leq \frac{1}{\sqrt{2}}|\tilde u-
\hat u|_{C^0([0,T_1];L^2(0,L))}.
\end{array}
$$
Hence, we get that $\cJ$ is a contractive
map. By the Banach fixed point theorem, we
know that $\cJ$ has a unique fixed point
$v_1$, which is the solution to the
following equation
\begin{equation}\label{12.19-eq2}
\left\{
\begin{array}{ll}\ds
v_{1,t} + v_{1,xxx} + v_{1,x} = \tilde h
&\mbox{ in }
[0,T_1]\times (0,L),\\
\ns\ds v_1(t,0)=v_1(t,L)=0 &\mbox{ on } [0,T_1],\\
\ns\ds v_{1,x}(t,L) = F(v_1(t,\cd)) &\mbox{ on } [0,T_1],\\
\ns\ds v_1(0,x)=v^0(x) &\mbox{ in } (0,L).
\end{array}
\right.
\end{equation}
Using \eqref{est1}, \eqref{defT1} and
\eqref{12.19-eq2}, we find that
$$
\begin{array}{ll}\ds
\q|v_1|_{C^0([0,T_1];L^2(0,L))}^2
\\
\ns\ds \leq C_1\(|v^0|_{L^2(0,L)}^2 +
\int_0^{T_1} \big|F(v(t,\cd))\big|^2 dt +
|\tilde h|_{L^1(0,T_1;L^2(0,L))}^2\)
\\
\ns \ds \leq C_1\(|v^0|_{L^2(0,L)}^2 + T_1
||F||^2 |v_1|_{C^0([0,T_1];L^2(0,L))}^2 +
|\tilde h|_{L^1(0,T_1;L^2(0,L))}^2\)
\\
\ns \ds \leq C_1\big(|v^0|_{L^2(0,L)}^2 +
|\tilde h|_{L^1(0,T_1;L^2(0,L))}^2\big) +
\frac{1}{2}
|v_1|_{C^0([0,T_1];L^2(0,L))}^2,
\end{array}
$$
which implies that
\begin{equation}\label{12.19-eq3}
|v_1|_{C^0([0,T_1];L^2(0,L))}^2 \leq
2C_1\big(|v^0|_{L^2(0,L)}^2 + |\tilde
h|_{L^1(0,T_1;L^2(0,L))}^2\big).
\end{equation}
From \eqref{est1}, \eqref{12.19-eq2} and
\eqref{12.19-eq3}, we obtain that
\begin{equation}\label{estim-v1}
\begin{array}{ll}\ds
\q|v_1|_{C^0([0,T_1];L^2(0,L))}^2 +
|v_1|^2_{L^2(0,T_1;H^1_0(0,L))} \\
\ns\ds\leq C_1\( |v^0|_{L^2(0,L)}^2 +
\int_0^{T_1} \big|F(v(t,\cd))\big|^2dt +
|\tilde h|_{L^1(0,T_1;L^2(0,L))}^2\)
\\
\ns \ds \leq C_1\big( |v^0|_{L^2(0,L)}^2 +
|\tilde h|_{L^1(0,T_1;L^2(0,L))}^2\big) +
C_1T_1 ||F||^2
|v_1|^2_{C^0([0,T_1];L^2(0,L))}
\\
\ns\ds \leq C_1\big( |v^0|_{L^2(0,L)}^2 +
|\tilde h|_{L^1(0,T_1;L^2(0,L))}^2\big) +
\frac{1}{2}|v_1|^2_{C^0([0,T_1];L^2(0,L))}
\\
\ns\ds \leq 2C_1\big( |v^0|_{L^2(0,L)}^2 +
|\tilde h|_{L^1(0,T_1;L^2(0,L))}^2\big).
\end{array}
\end{equation}
By a similar argument, we can prove that
the following equation
\begin{equation}\label{12.19-eq2.1}
\left\{
\begin{array}{ll}\ds
v_{2,t} + v_{2,xxx} + v_{2,x} = \tilde h
&\mbox{ in }
[T_1,2T_1]\times (0,L),\\
\ns\ds v_2(t,0)=v_2(t,L)=0 &\mbox{ on } [T_1,2T_1],\\
\ns\ds v_{2,x}(t,L) = F(v_2(t,\cd)) &\mbox{ on } [T_1,2T_1],\\
\ns\ds v_2(T_1,x))= v_1(T_1,x)&\mbox{ in }
(0,L).
\end{array}
\right.
\end{equation}
admits a unique solution in
$C^0([T_1,2T_1];L^2(0,L))\cap
L^2(T_1,2T_1;H^1_0(0,L))$. Furthermore,
this solutions satisfies that
\begin{equation}\label{estim-v2}
\begin{array}{ll}\ds
\q|v_2|_{C^0([T_1,2T_1];L^2(0,L))}^2 +
|v_2|^2_{L^2(T_1,2T_1;H^1_0(0,L))}
\\
\ns\ds \leq 2C_1\big(
|v_1(T_1)|_{L^2(0,L)}^2 + |\tilde
h|_{L^1(T_1,2T_1;L^2(0,L))}^2\big).
\end{array}
\end{equation}
By induction, we know that for an integer
$n\geq 2$, the following equation admits a
unique solution.
\begin{equation}\label{eqvn-A}
\left\{
\begin{array}{ll}\ds
v_{n,t} + v_{n,xxx} + v_{n,x} = \tilde h
&\mbox{ in }
[(n-1)T_1,nT_1]\times (0,L),\\
\ns\ds v_n(t,0)=v_n(t,L)=0 &\mbox{ on } [(n-1)T_1,nT_1],\\
\ns\ds v_{n,x}(t,L) = F(v_n(t,\cd)) &\mbox{ on } [(n-1)T_1,nT_1],\\
\ns\ds v_n((n-1)T_1,x)=
v_{n-1}((n-1)T_1,x)&\mbox{ in } (0,L).
\end{array}
\right.
\end{equation}
Furthermore, one has
\begin{equation}\label{estim-vn}
\begin{array}{ll}\ds
\q|v_n|_{C^0([(n-1)T_1,nT_1];L^2(0,L))}^2 +
|v_n|^2_{L^2((n-1)T_1,nT_1;H^1_0(0,L))}
\\
\ns\ds \leq 2C_1\big(
|v_{n-1}((n-1)T_1)|_{L^2(0,L)}^2 + |\tilde
h|_{L^1((n-1)T_1,nT_1;L^2(0,L))}^2\big).
\end{array}
\end{equation}
Let $n_0\in \dbZ^+$ be such that
$(n_0-1)T_1 \leq T < n_0T_1$. Let
$$
v(t,x)\=v_n(t,x) \mbox{ for } t\in
[(n-1)T_1,nT_1)\cap [0,T],\q
n=1,2,\cds,n_0, \, x\in (0,L).
$$
Then $v$ is the solution to
\begin{equation}\label{eqv-A}
\left\{
\begin{array}{ll}\ds
v_{t} + v_{xxx} + v_{x} = \tilde h &\mbox{
in }
[0,T]\times (0,L),\\
\ns\ds v(t,0)=v(t,L)=0 &\mbox{ on }
[0,T],\\
\ns\ds v_{x}(t,L) = F(v(t,y)) &\mbox{ on } [0,T],\\
\ns\ds v(0,x)= v^0(x)&\mbox{ in } (0,L).
\end{array}
\right.
\end{equation}
Moreover,  there is a constant $C=C(T)>0$
(independent of $v^0\in L^2(0,L)$ and of
$\tilde h\in L^1(0,T;L^2(0,L))$) such that
\begin{equation}\label{estim-lin-kdv}
|v|_{C^0([0,T];L^2(0,L))}^2 +
|v|^2_{L^2(0,T;H^1_0(0,L))} \leq C(T)\big(
|v^0|_{L^2(0,L)}^2 + |\tilde
h|_{L^1(0,T;L^2(0,L))}^2\big).
\end{equation}

 Let $\kappa>0$, depending only on $L>0$, such that
\begin{equation}\label{Sobolev}
  |w|_{L^\infty}\leq \kappa |w|_{H_0^1(0,L)},\q \forall w\in H_0^1(0,L).
\end{equation}
 From now on, we assume that $v^0\in L^2(0,L)$ satisfies
\begin{equation}\label{conditionu0}
  |v^0|_{L^2(0,L)}^2 \leq
\frac{5}{36C(T)^2\kappa^2}.
\end{equation}
 Let
$$
\begin{array}{ll}\ds
\cB\= \Big\{u \in C^0\!([0,T];
L^2(0,L))\cap L^2(0,T;H^1_0(0,L)):
\\
\ns\ds \qq\q |u|^2_{C^0([0,T];L^2(0,L))}+
|u|^2_{L^2(0,T;H^1_0(0,L))} \leq
\frac{1}{6C(T)\kappa^2} \Big\}.
\end{array}
$$
Then $\cB$ is a nonempty closed subset of
the Banach space $C^0\!([0,T];
L^2(0,L))\cap L^2(0,T;H^1_0(0,L))$. Let us
define a map $\cK$ from $\cB$ to
$C^0([0,T];L^2(0,L))\cap
L^2(0,T;H^1_0(0,L))$ as follows:
$$
\cK(z) = v, \mbox{ where $v$ is the
solution to \eqref{eqv-A} with } \tilde h =
zz_x.
$$
From Lemma \ref{lm4} and the above
well-posedness result for \eqref{eqv-A}, we
know that $\cK$ is well-defined. From
\eqref{estim-lin-kdv} and
\eqref{conditionu0}, we have that
$$
\begin{array}{ll}\ds
\q|\cK(z)|^2_{C^0([0,T];L^2(0,L))} +
|\cK(z)|^2_{L^2(0,T;H^1_0(0,L))}
\\
\ns\ds \leq C(T)\big(|v^0|_{L^2(0,L)}^2 +
|zz_x|_{L^1(0,T;L^2(0,L))}^2\big)
\\
\ns\ds \leq C(T)\big(|v^0|_{L^2(0,L)}^2 +
|z|^2_{L^2(0,T;L^\infty(0,L))}|z_x|^2_{L^2(0,T;L^2(0,L))}\big)
\\
\ns\ds \leq C(T)\big( |v^0|_{L^2(0,L)}^2 +
\kappa^2 |z|^4_{L^2(0,T;H^1_0(0,L))}\big)
\\
\ns\ds \leq \frac{1}{6C(T)\kappa^2}
\end{array}
$$
and
$$
\begin{array}{ll}\ds
\q|\cK(z_1)-\cK(z_2)|^2_{C^0([0,T];L^2(0,L))}
+
|\cK(z_1)-\cK(z_2)|^2_{L^2(0,T;H^1_0(0,L))}
\\
\ns\ds \leq C(T)
|z_1z_{1,x}-z_2z_{2,x}|^2_{L^1(0,T;L^2(0,L))}=
C(T) |z_1z_{1,x}-z_1z_{2,x} +
z_1z_{2,x}-z_2z_{2,x}
|^2_{L^1(0,T;L^2(0,L))}
\\
\ns\ds \leq 2C(T)\kappa^2 \big(
|z_1|^2_{L^2(0,T;H_0^1(0,L))}|z_1-z_2|^2_{L^2(0,T;H_0^1(0,L))}
+
|z_1-z_2|^2_{L^2(0,T;H_0^1(0,L))}|z_2|^2_{L^2(0,T;H_0^1(0,L))}
\big)
\\
\ns\ds \leq
2C(T)\kappa^2\big(|z_1|^2_{L^2(0,T;H_0^1(0,L))}+|z_2|^2_{L^2(0,T;H^1_0(0,L))}\big)
|z_1-z_2|^2_{L^2(0,T;H^1_0(0,L))} \\
\ns\ds \leq
\frac{2}{3}|z_1-z_2|^2_{L^2(0,T;H^1_0(0,L))}.
\end{array}
$$
Therefore, we know that $\cK$ is from $\cB$
to $\cB$ and is contractive. Then, by the
Banach fixed point theorem, there is a
(unique) fixed point $v$, which is the
solution to \eqref{csystem1} with initial
data $v(0,\cdot)=v^0(\cdot)$.

\begin{remark}
\label{largeinitialdata} Adapting the proof
of \cite[Section 3.1]{MVZ}, one can also
prove that, for every $v^0\in L^2(0,L)$,
there exist $T>0$ and a solution  $v\in
C^0([0,T]; L^2(0,L))\cap L^2(0,T;H^1(0,L))$
of \eqref{csystem1} with $f(t)=F(v(t,\cd))$
satisfying the initial condition
$v(0,\cdot)=v^0(\cdot)$.
\end{remark}


\section*{Acknowledgement} Jean-Michel
Coron was supported by ERC advanced grant
266907 (CPDENL) of the 7th Research
Framework Programme (FP7).

Qi L\"{u} was supported by ERC advanced
grant 266907 (CPDENL) of the 7th Research
Framework Programme (FP7), the NSF of China
under grant 11101070 and the Fundamental
Research Funds for the Central Universities
in China under grants ZYGX2012J115. Qi
L\"{u} would like to thank the Laboratoire
Jacques-Louis Lions (LJLL) at
Universit\'{e} Pierre et Marie Curie for
its hospitality. This work was carried out
while he was visiting the LJLL.

\bibliographystyle{plain}
\bibliography{kernel}

\end{document}